\documentclass[11pt]{article}
\usepackage{amssymb,amsmath}% font used for R in Real numbers
\usepackage{mathrsfs}
\usepackage{latexsym}
\usepackage{color}
\usepackage{amssymb}
\usepackage{amsmath}
\usepackage{amsthm}
\usepackage{amsmath,amsfonts,amsthm,amssymb,amscd,color,xcolor,mathrsfs,verbatim,microtype}

\usepackage{graphicx,eurosym}
\usepackage{hyperref}
\usepackage[applemac]{inputenc}
\usepackage[cyr]{aeguill}
\colorlet{darkblue}{blue!50!black}
\hypersetup{
    colorlinks,%
    citecolor=blue,%
    filecolor=red,%
    linkcolor=darkblue,%
    urlcolor=blue,%
    pdfnewwindow=true,%
    pdfstartview={FitH}
}

\allowdisplaybreaks

\catcode`!=11
\let\!int\int \def\int{\displaystyle\!int}
\let\!lim\lim \def\lim{\displaystyle\!lim}
\let\!sum\sum \def\sum{\displaystyle\!sum}
\let\!sup\sup \def\sup{\displaystyle\!sup}
\let\!inf\inf \def\inf{\displaystyle\!inf}
\let\!cap\cap \def\cap{\displaystyle\!cap}
\let\!max\max \def\max{\displaystyle\!max}
\let\!min\min \def\min{\displaystyle\!min}
\catcode`!=12

\newtheorem{theorem}{\bf Theorem}[section]
\newtheorem{lemma}{\bf Lemma}[section]
\newtheorem{definition}{\bf Definition}[section]
\newtheorem{proposition}{\bf Proposition}[section]

\newtheorem{remark}{\bf Remark}[section]

\renewenvironment{abstract}{%
  \small
  \begin{center}%
    {\bfseries \abstractname\vspace{-.5em}\vspace{0pt}}%
  \end{center}%
  \list{}{%
    \setlength{\leftmargin}{2.3em}%
    \setlength{\rightmargin}{\leftmargin}%
  }%
  \item\relax
}{%
  \endlist
}

\catcode`@=11
\@addtoreset{equation}{section}
\catcode`@=12

\begin{document}
\title{Polynomial mixing for the white-forced wave equation on the whole line}
\author{Peng Gao
\\[2mm]
\small School of Mathematics and Statistics, and Center for Mathematics
\\
\small and Interdisciplinary Sciences, Northeast Normal University,
\\
\small Changchun 130024,  P. R. China
\\[2mm]
\small Email: gaopengjilindaxue@126.com }
\date{\today}
\maketitle

\vbox to -13truemm{}

\begin{abstract}
Our goal in this paper is to investigate the ergodicity of the white-forced wave equation on the whole line. Assuming that sufficiently many directions of the phase space are stochastically forced, we prove the uniqueness of the stationary measure and polynomial mixing in the dual-Lipschitz metric. The difficulties in our proof are twofold. On the one hand, in contrast to stochastic parabolic equations, the stochastic wave equation lacks both a smoothing effect and a strongly dissipative mechanism. On the other hand, the unboundedness of the whole line leads to a loss of compactness compared with the bounded domain case. To overcome these obstacles, our proof relies on four key ingredients: a new criterion for polynomial mixing established in \cite{Gao26JDE}, a novel weighted Foia\c{s}-Prodi estimate for the wave equation on the whole line, weighted energy estimates for the stochastic wave equation, and the irreducibility of the stochastic wave equation.
\\[6pt]
\end{abstract}

\tableofcontents
\setcounter{section}{0}

\section{Introduction}
\subsection{Motivation}
Wave equation
$$\ddot{u}-\Delta u+f(u)=0$$
is an important class of partial differential equations, describing a great variety of wave phenomena, it appears in the study of several problems of mathematical physics. For example, this equation arises in general relativity, nonlinear optics (e.g., the instability phenomena such as self-focusing), plasma physics, fluid mechanics, radiation theory or spin waves. Since the pioneering
work of Lions and Strauss \cite{Lions1965}, the deterministic wave equation has been extensively investigated, the developments in various directions are well documented in \cite{Har1987,Che2008,Che2017,Chu,Cot2021,Fei1995,Zua1990,Zua1991,Zua1993,Zua2003} and the references therein.
Stochastic wave equation arises as a mathematical
model to describe nonlinear vibration or wave propagation in a randomly excited
continuous medium, such as atmosphere, oceans, sonic booms, traffic flows, optic
devices and quantum fields, when random fluctuations are taken into account.
Over the past two decades, stochastic wave equation has become an active research area which attracted much attention, it has been studied by many
authors, see \cite{Barbu1,Barbu2,Chow1,Chow2,Chow2015,Kim2004} and the references therein.
\par
In this paper, we consider the long time behavior for the following white-forced wave equation on the whole line
\begin{eqnarray}\label{WE}
 \begin{array}{l}
 \left\{
 \begin{array}{llll}
\ddot{u}+Au+\gamma \dot{u}+f(u)=h+\eta
 \\u(x,0)=u_{0}(x)
 \\\dot{u}(x,0)=u_{1}(x)
 \end{array}
 \right.
 \end{array}
 \begin{array}{lll}
 {\rm{in}}~\mathbb{R}\times(0,+\infty),\\
 {\rm{in}}~\mathbb{R},
 \\{\rm{in}}~\mathbb{R},
\end{array}
\end{eqnarray}
where $\dot{u}=\frac{du}{dt}, \ddot{u}=\frac{d^{2}u}{dt^{2}}, Au=(-\partial_{xx}+1)u, f(u)=|u|^{2m}u~(0<m<1),$ $\gamma>0$ is a constant, $h=h(x),$ $\eta=\eta(x,t)$ is a white noise of the form
\begin{equation}\label{37}
\begin{split}
\eta(x,t):=\frac{\partial}{\partial t}W(t),~~W(t):=\sum\limits_{i=1}^{\infty}b_{i}\beta_{i}(t)e_{i}(x),
\end{split}
\end{equation}
and $b_{i}\in \mathbb{R}, \{\beta_{i}\}_{i\geq1}$ is a sequence of
independent real-valued standard Brownian motions defined on a filtered probability space
$(\Omega,\mathcal{F},\mathcal{F}_{t},\mathbb{P})$ satisfying the usual conditions, $\{e_{i}\}_{i\geq 1}$ is an orthonormal basis in $L^{2}(\mathbb{R}).$
Let $H := L^2(\mathbb{R})$ and $\mathcal{H} := H^1(\mathbb{R}) \times L^2(\mathbb{R})$. We equip $\mathcal{H}$ with its usual norm
\[
\|y\|_{\mathcal{H}}^2 := \|y_1\|_{1}^2 + \|y_2\|^2, \qquad y = (y_1, y_2) \in \mathcal{H}.
\]
However, for the purposes of this paper, it is convenient to introduce the following equivalent norm:
\[
|y|_{\mathcal{H}}^2 := \|y_1\|_{1}^2 + \|y_2 + \alpha y_1\|^2, \qquad y = (y_1, y_2) \in \mathcal{H},
\]
where $\alpha > 0$ is a sufficiently small constant. It is straightforward to verify that the norms $\|\cdot\|_{\mathcal{H}}$ and $|\cdot|_{\mathcal{H}}$ are equivalent, provided that $\alpha > 0$ is taken to be sufficiently small.
\par
Let $y_0=[u_0,u_1]$, $y^u=[u,\dot{u}]$, for notational simplicity, we occasionally write $y$ instead of $y^u$.
Setting
$\mathbb{F}(y)=[0,-f(u)+h], \mathbb{G}=[0,1]$ and
$
\mathbb{A}=\left(
  \begin{array}{cc}
    0 & 1 \\
    -A & -\gamma \\
  \end{array}
\right),
$
we can rewrite \eqref{WE} as
\begin{eqnarray}\label{WEY}
 \begin{array}{l}
 \left\{
 \begin{array}{llll}
dy=[\mathbb{A}y+\mathbb{F}(y)]dt+\mathbb{G}dW
 \\y(0)=y_{0}
 \end{array}
 \right.
 \end{array}
 \begin{array}{lll}
 {\rm{in}}~\mathbb{R}\times(0,+\infty),\\
 {\rm{in}}~\mathbb{R}.
\end{array}
\end{eqnarray}
It is well known (see \cite{Pazy}) that $\mathbb{A}$ is the infinitesimal generator of a strongly
continuous group in $\mathcal{H}$ given by $\{U(t)\}_{t\in \mathbb{R}}$. With the help of the group $\{U(t)\}_{t\in \mathbb{R}}$,
we can give the following definition of the mild solution to system \eqref{WE}.
\begin{definition}\label{mild}
Given $(\Omega,\mathcal{F},\mathcal{F}_{t},\mathbb{P})$ and an $\mathcal{F}_{0}$-measurable real $\mathcal{H}-$valued initial condition $y_{0}$, a process $y$ is called a \textit{mild solution} of (\ref{WE})
on $[0,T],$ if for almost each $\omega\in \Omega$ and $t\in[0,T],$ $y$ satisfies the following
It\^{o} integral form
\begin{equation*}
\begin{split}
y(t)=U(t)y_{0}+\int_{0}^{t}U(t-s)\mathbb{F}(y(s))ds+\int_{0}^{t}U(t-s)\mathbb{G}dW(s).
\end{split}
\end{equation*}
\end{definition}
\par
We can show that the stochastic wave equation \eqref{WE} is globally well-posed, its proof can be found in Appendix.
The definition of $\mathcal{E}_u(t)$ is given in next section.
\begin{proposition}\label{WP}
Let $h\in H$ and $y_0$ be an $\mathcal{H}-$valued random variable that is independent of $W$ and satisfies $\mathbb{E}\mathcal{E}_u(0)<+\infty$. Then system \eqref{WE} possesses a unique solution in the sense of Definition \ref{mild}. In addition, we have the a priori estimate
\begin{equation}\label{45}
\begin{split}
\mathbb{E}\mathcal{E}_u(t)\leq \mathbb{E}\mathcal{E}_u(0)e^{-c t}+C,
\end{split}
\end{equation}
where $c,C>0$ are some constants.
\end{proposition}

\par
By the standard arguments as in \cite{Barbu1,Barbu2,Chow1,Chow2,DaZ1,M2014}, we can define a Markov process in $\mathcal{H}$. Let the transition semigroup $\mathcal{B}_t$ on $C_{b}(\mathcal{H})$ associated with the flow $t\mapsto y(t)=[u(t),\dot{u}(t)]$ be
$$(\mathcal{B}_tf)(y_{0}):=\mathbb{E}f(y(t)),~~f\in C_{b}(\mathcal{H}),$$
here $C_{b}(\mathcal{H})$ is the Banach space of all uniformly continuous and bounded mappings on $\mathcal{H}$.
Let $(y(t),\mathbb{P}_{y})$ be the Markov process in $\mathcal{H}$ corresponding to \eqref{WE},
then
$$P_{t}(y,A):=\mathbb{P}(y_{t}(y,\omega)\in A)$$
is the transition function. We introduce the following associated Markov operators
\begin{equation*}
\begin{split}
&\mathcal{B}_{t}: C_{b}(\mathcal{H})\rightarrow C_{b}(\mathcal{H}),~\mathcal{B}_{t}f(y):=\int_{\mathcal{H}}f(z)P_{t}(y,dz),~\forall f\in C_{b}(\mathcal{H}),\\
&\mathcal{B}_{t}^{*}: \mathcal{P}(\mathcal{H})\rightarrow \mathcal{P}(\mathcal{H}),~\mathcal{B}_{t}^{*}\lambda(A):=\int_{\mathcal{H}}P_{t}(y,A)\lambda(dy),~\forall \lambda\in \mathcal{P}(\mathcal{H}).
\end{split}
\end{equation*}
\par
Motivated from both physical and mathematical standpoints, an important mathematical question arises:
\par
~~
\par
\rm
\textbf{What is the asymptotic behavior of $u(t)$ as $t\rightarrow+\infty$~?}
\rm
\par
~~
\par
The asymptotic behavior of solutions to the wave equation is a classical topic, and we briefly recall some relevant contributions. Energy decay and stabilization have been investigated in \cite{Zua1990,Zua1991,Zua2003,Zua2006,Cot2021}; attractors for the deterministic wave equation have been studied in \cite{Zua1993,Fei1995,Chu,Zel2009}; and averaging principles have been established in \cite{Che2008,Che2017}. In recent years, the long-time dynamics of stochastic wave equations have attracted increasing interest; see, for example, \cite{Chow1,Chow2,Chow2015,PMAMO,att2015,att2011} and the references therein. Among various approaches, the mixing property has proven to be an effective tool for describing the long-time behavior of stochastic partial differential equations (SPDEs). We refer the reader to \cite{Deb1,HM06,HM08,HM11-1} and the monographs \cite{DaZ2,KS12} for a detailed account. In this paper, we establish polynomial mixing for the solution of system \eqref{WE}.

\subsection{Main result}
\par
Now, we are in a position to present the main result in this paper, this can answer the above question. We introduce the following assumption
\begin{equation*}
\begin{array}{l}
\textbf{(A)}
\left\{
\begin{array}{lll}
&h, e_i\in H^2,\varphi h, \varphi e_i\in H, i\geq 1,~\sum\limits_{i=1}^{\infty}|(h,e_{i})|\|e_{i}\|_{2}<+\infty,\\
&\mathcal{B}_{1}:=\sum\limits_{i=1}^{\infty}b_{i}^{2}<+\infty,\\
&\mathcal{B}_{2}:=\sum\limits_{i=1}^{\infty}b_{i}^{2}\|\varphi e_{i}\|^{2}<+\infty,\\
&\mathcal{B}_{3}:=\sum\limits_{i=1}^{\infty}|b_{i}|\|e_{i}\|_{2}<+\infty,
\end{array}
\right.
\end{array}
\end{equation*}
where $\varphi(x):=\ln (x^{2}+2)$.
We state the result for polynomial mixing of system \eqref{WE}.
\begin{theorem}\label{MT}
Let the assumption $\textbf{(A)}$ hold. Then, there exists an integer $N\geq1$ such that if
\begin{equation}\label{41}
\begin{split}
b_{i}\neq0,~i=1,2,\cdots,N,
\end{split}
\end{equation}
then there exists a unique stationary measure $\mu\in \mathcal{P}(\mathcal{H})$ for \eqref{WE}. Moreover, for any $p>1$, there exist an integer $N=N_p\geq1$ (for which \eqref{41} holds)  and a positive constant $C_p$ such that for any $y_{0}=[u_0,u_1]\in \mathcal{H},$ the solution $y$ of \eqref{WE} satisfies
\begin{equation}\label{0}
\begin{split}
|\mathbb{E}f(y(t))-\int_{\mathcal{H}}f(y)\mu(dy)|\leq C_p(t+1)^{-p}\|f\|_{Lip}(1+\|u_0\|_1^2+\|u_1\|^2+\|u_0\|^{2m+2}_{L^{2m+2}}),~t\geq 0,
\end{split}
\end{equation}
for any bounded Lipschitz-continuous function $f:\mathcal{H}\to \mathbb{R}$.

\end{theorem}

\begin{remark}
The estimate \eqref{0} is equivalent to the following estimate
\begin{equation*}
\begin{split}
\|P_{t}(y_{0},\cdot)-\mu\|_{L(\mathcal{H})}^{*}\leq C_p(t+1)^{-p}(1+\|u_0\|_1^2+\|u_1\|^2+\|u_0\|^{2m+2}_{L^{2m+2}}),~t\geq 0,
\end{split}
\end{equation*}
and also to
\begin{equation*}
\begin{split}
\|\mathcal{B}_{t}^{*}\lambda-\mu\|_{L(\mathcal{H})}^{*}\leq C_p(t+1)^{-p}(1+\int_{\mathcal{H}}(\|u\|_1^2+\|v\|^2+\|u\|^{2m+2}_{L^{2m+2}})\lambda(d(u,v))),~t\geq 0,
\end{split}
\end{equation*}
for any $\lambda\in \mathcal{P}(\mathcal{H}).$ We refer the reader to the end of Section~1 for the definitions of $\|\cdot\|^*_{L(\mathcal{H})},$ $\|\cdot\|_{Lip}, \|\cdot\|_s$ and $\|\cdot\|_{L^p}$.
\end{remark}
\par~~
\par
In the last thirty years, there was a substantial progress in ergodicity for partial differential equations with random forcing, we refer the readers
to \cite{Deb1,HM06,HM08,HM11-1} and the books \cite{DaZ2,KS12} for a detailed discussion of this direction. In the beginning, people mainly focus on mixing for parabolic type SPDEs, such as Navier-Stokes system, Ginzburg-Landau equation, reaction-diffusion equation and so on. Most of the existing literature has been mainly concerned on the case of exponential mixing for parabolic type SPDEs, we refer readers to \cite{KS12,KNS1,Shi08,HM08,HM11-1,FG15} and the references therein.
It is important to point out that in the above works, the exponential convergent rate is valid thanks to the advantage of strong dissipation. Generally speaking, it is difficult to be achieved for weakly dissipative SPDEs.
\par
In recent years, mixing problem for weakly dissipative SPDEs has received more and more attention, this kind of SPDEs does not possess a strong dissipative mechanism, it seems to be hard to study the mixing problem for this kind of SPDEs. Stochastic wave equation is an important kind of weakly dissipative SPDEs, it lacks the smoothing effect and strongly dissipative mechanism compared to parabolic type SPDEs, these reasons make the uniqueness of stationary
measure and its ergodicity for stochastic wave equation be much more delicate. Up to now, mixing for stochastic wave equation is studied less, here we recall them. \cite{Barbu1,Barbu2} established the uniqueness of the invariant measure for the white-forced wave equation with the help of Kolmogorov operator, \cite{Chow2} investigated the uniqueness of the invariant measure for the wave equation by energy dissipativity method, \cite{PMAMO} investigated polynomial mixing for white-forced wave equation with the help of asymptotic coupling method, \cite{Seong2024,For2024} investigated the uniqueness of the invariant measure for wave equation with space-time white noise and \cite{M2014} established exponential mixing for the white-forced damped nonlinear wave equation by virtue of coupling method. All the above results are established on a \textit{bounded domain}.
To the best of our knowledge, there are rarely studies on the mixing for stochastic wave equation on an \textit{unbounded domain}.
\par
Up to now, there is less work on the mixing of SPDEs on the whole line. Compared to bounded
domains, the whole line will lead to the lack of compactness.
\cite{N2} established a beautiful framework and succeeded in establishing exponential mixing for the white-forced complex Ginzburg-Landau equation by virtue of coupling method from \cite{KS12,Shi08}.
\cite{Gao26JDE} investigated polynomial mixing for the white-forced hyperviscous Burgers equation on the whole line with the help of a new coupling criterion (see \cite[Theorem 1.2 and Remark 1.1]{Gao26JDE}). Motivated by \cite{N2,Gao26JDE}, in this paper, we establish polynomial mixing of the white-forced wave equation on the whole line. This paper is an attempt for mixing of weakly dissipative SPDEs on the whole line.

We employ the coupling method to prove Theorem \ref{MT}. The main difficulties in the proof are twofold. First, unlike stochastic parabolic equations, the stochastic wave equation possesses neither a smoothing effect nor a strong dissipative mechanism. Second, the unboundedness of the entire real line causes a loss of compactness compared with the bounded domain setting. To overcome these challenges, our argument relies on four essential ingredients:
\par
$\bullet$ A new criterion for polynomial mixing established in \cite{Gao26JDE} (see Theorem \ref{Th1}). This criterion weakens the conditions of the previous exponential mixing criterion, thereby making it applicable to a wider class of random dynamical systems. It also serves as a tool for establishing mixing rates for weakly dissipative SPDEs.
\par
$\bullet$ A novel weighted Foia\c{s}-Prodi estimate for the stochastic wave equation on the whole line (see Theorem \ref{FP}). The Foia\c{s}--Prodi estimate plays a crucial role in characterizing the asymptotic behavior of solutions to SPDEs emanating from different initial data. It also serves as an important tool for verifying the conditions of the abstract criterion (see Theorem \ref{Th1}). Compared with the bounded domain case, we introduce a weight function to overcome the lack of compactness; this is motivated by \cite{N2}.
\par
$\bullet$ Weighted energy estimates for the stochastic wave equation developed in Section 3. In order to apply the Foia\c{s}-Prodi estimate to establish the polynomial mixing theory via the coupling method, we need to establish suitable weighted energy estimates for the wave equation on the whole line. Unlike the parabolic case, the corresponding energy estimates for wave equations are considerably more involved and require more sophisticated techniques.
\par
$\bullet$ Irreducibility of the stochastic wave equation (see Proposition \ref{pro20}). Irreducibility is another critical ingredient for establishing mixing. By combining the continuity of the solution operator with the support property of the noise, we prove that the Markov transition probabilities of the stochastic wave equation on the whole line enjoy a certain irreducibility property.

We believe that the strategy and framework developed in this paper can be applied to more complicated weakly dissipative SPDEs, and that the mixing problem can be treated via the corresponding Foia\c{s}--Prodi estimate, weighted energy estimates and irreducibility.
\par~~
\par
The rest of the paper is organized as follows. In Section 2, we introduce the coupling method and an abstract criterion for the proof of Theorem \ref{MT} and establish a new Foia\c{s}-Prodi estimate
for wave equation on the whole line.
We establish some estimates for solution to \eqref{WE} in Section 3.
In Section 4, we prove the irreducibility for \eqref{WE}. Finally, we prove Theorem \ref{MT}.
Some auxiliary results are proved in Appendix.
\par~~
\par~~
\par
\textbf{Notation}
\par
Throughout the paper, $c$ and $C$ denote generic positive constants that may change
from line to line. We adopt the convention that constants are written with their parameters in parentheses or as subscripts; for example, $C(p,T)$ depends on $p$ and $T$, and $C_p$ depends on $p$.
\par
Let $X$ be a Polish space with a metric $d_{X}(u,v)$, the Borel $\sigma$-algebra on $X$ is denoted by $\mathcal{B}(X)$ and the set of Borel
probability measures by $\mathcal{P}(X).$ We denote by $B_X(x,r)$ the open ball in $X$ with radius $r$ centered at $x\in X$, and by $\bar{B}_X(x,r)$ its closure. $C_{b}(X)$ is the space of continuous functions $f:X\rightarrow \mathbb{C}$ endowed with the
norm $\|f\|_{\infty}=\sup_{u\in X}|f(u)|.$
We write $C(X)$ when $X$ is compact. $L_{b}(X)$ is the space of functions $f\in C_{b}(X)$ such that
$$\|f\|_{Lip}=\|f\|_{\infty}+\sup\limits_{u\neq v}\frac{|f(u)-f(v)|}{d_{X}(u,v)}<\infty.$$
The dual-Lipschitz metric on $\mathcal{P}(X)$ is defined by
$$\|\mu_{1}-\mu_{2}\|^{*}_{L(X)}=\sup_{\|f\|_{L(X)}\leq 1}|\langle f,\mu_{1}\rangle-\langle f,\mu_{2}\rangle|,~~\mu_{1},\mu_{2}\in \mathcal{P}(X),$$
where $\langle f,\mu \rangle=\int_{X}f(u)\mu(du).$ The total variation metric on $\mathcal{P}(X)$ is defined by
\begin{equation*}
\|\mu_{1}-\mu_{2}\|_{var}=\frac12\sup_{f\in C_b(X),\|f\|_{\infty}\leq1}|\langle f,\mu_{1}\rangle-\langle f,\mu_{2}\rangle|,~~\mu_{1},\mu_{2}\in \mathcal{P}(X).
\end{equation*}
\par
Let $H^{s}:=H^{s}(\mathbb{R})$ denote the Sobolev space of real-valued functions on $\mathbb{R}$, endowed with the standard Sobolev norm $\|\cdot\|_{s}$ (see \cite{Adams}). In the particular case $s=0$, we write $H:=L^{2}(\mathbb{R})$, denote its norm simply by $\|\cdot\|$, and its inner product by $(\cdot,\cdot)$. Moreover, for $1\le p\le \infty$, let $L^{p}:=L^{p}(\mathbb{R})$ be the usual Lebesgue space of $p$-th power integrable functions on $\mathbb{R}$, equipped with the canonical norm $\|\cdot\|_{L^{p}}$ (see \cite{Adams}).
\par
For $s\geq0,$ we define $\mathcal{G}^{s}:=\{\zeta\in C([0,T];H^{s})~|~\zeta(0)=0\}.$
\par
For $q>1$, let $Q_q$ be the function defined by $Q_q(r):=r^{1-q}$ for $r>0$.

\section{Coupling method}
\subsection{An abstract criterion}
\par
Let $X$ be a separable Banach space with a norm $\|\cdot\|$. Let $(u_{t}, \mathbb{P}_{u})$ be a Feller family of Markov processes in $X,$ $P_{t}(u,A):=\mathbb{P}(u_{t}(u,\omega)\in A)$ is
the transition function. We introduce the following associated Markov operators
\begin{equation*}
\begin{split}
&\mathcal{B}_{t}: C_{b}(X)\rightarrow C_{b}(X),~\mathcal{B}_{t}f(u):=\int_{X}f(v)P_{t}(u,dv),~\forall f\in C_{b}(X),\\
&\mathcal{B}_{t}^{*}: \mathcal{P}(X)\rightarrow \mathcal{P}(X),~\mathcal{B}_{t}^{*}\lambda(A):=\int_{X}P_{t}(u,A)\lambda(du),~\forall \lambda\in \mathcal{P}(X).
\end{split}
\end{equation*}

\begin{definition}(See \cite{KS12,Shi08})
Let $(\mathbf{u}_{t},\mathbb{P}_{\mathbf{u}})$ be a Markov family in $X\times X.$
$(\mathbf{u}_{t},\mathbb{P}_{\mathbf{u}})$ is called an \textbf{extension} of $(u_{t}, \mathbb{P}_{u})$ if for any $\mathbf{u}=(u,u^{\prime})\in X\times X$ the laws under $\mathbb{P}_{\mathbf{u}}$ of processes $\{\Pi_{1}\mathbf{u}_{t}\}_{t\geq0}$ and $\{\Pi_{2}\mathbf{u}_{t}\}_{t\geq0}$ coincide with
those of $\{u_{t}\}_{t\geq0}$ under $\mathbb{P}_{u}$ and $\mathbb{P}_{u^{\prime}}$ respectively, where $\Pi_{1}$ and $\Pi_{2}$ denote the
projections from $X\times X$ to the first and second component.
\end{definition}

Now, we are in a position to present an abstract criterion for polynomial mixing from \cite{Gao26JDE}, its proof and exact form can be found in \cite[Theorem 1.2 and Remark 1.1]{Gao26JDE}.
This criterion is developed in a sufficiently general framework to be applied and prove polynomial convergence to equilibrium.
\begin{theorem}\label{Th1}\cite[Theorem 1.2 and Remark 1.1]{Gao26JDE}
Let $(u_{t},\mathbb{P}_{u})$ be a Feller family of Markov processes in $X$.
Supposed that
\par
(1) There exists an extension $(\mathbf{u}_{t},\mathbb{P}_{\mathbf{u}})$, a stopping time $\sigma,$ a closed set $B\subset X,$ and an increasing
function $g(r)\geq1$ of the variable $r\geq0$ such that the following two properties
hold
\par
\textbf{Recurrence} There exists $p>1$ such that
$$\mathbb{E}_{\mathbf{u}} \tau_{\textbf{B}}^{p}\leq G(\mathbf{u})$$
for all $\mathbf{u}=(u,u^{\prime})\in \textbf{X}:=X\times X,$
where $\tau_{\textbf{B}}=\tau_{\textbf{B}}(\textbf{u},\omega)=\inf\{t\geq0:\textbf{u}_{t}(\textbf{u},\omega)\in \textbf{B}:=B\times B\}$ and $G(\mathbf{u}):=g(\|u\|)+g(\|u^{\prime}\|)$.
\par
\textbf{Polynomial squeezing} There exist positive constants $c_{1},c_{2},c_3,c_4$ and $1\leq q\leq p$
such that, for any $\mathbf{u}\in \mathbf{B},$ we have
\begin{equation*}
\begin{split}
&\|u_{t}-u_{t}^{\prime}\|\leq c_{1}\|u-u^{\prime}\|(t+1)^{-p},~{\rm{for}}~0\leq t\leq \sigma,\\
&\mathbb{P}_{\mathbf{u}}(\sigma=\infty)\geq c_{2},\\
&\mathbb{E}_{\mathbf{u}}(\mathbb{I}_{\sigma<\infty}\sigma^{p})\leq c_3,\\
&\mathbb{E}_{\mathbf{u}}(\mathbb{I}_{\sigma<\infty}G(\mathbf{u}_{\sigma})^{q})\leq c_4.
\end{split}
\end{equation*}

(2) There is an increasing function $\tilde{g}(r)\geq1$ such that
\begin{equation*}
\begin{split}
\mathbb{E}_{u}g(\|u_{t}\|)\leq \tilde{g}(\|u\|),~{\rm{for}}~t\geq 0,~u\in X.
\end{split}
\end{equation*}

Then the family $(u_{t},\mathbb{P}_{u})$ has a unique stationary measure $\mu\in \mathcal{P}(X),$ and it holds for any $p_0\le p$ that
\begin{equation*}
\begin{split}
\|P_{t}(u,\cdot)-\mu\|^{*}_{L}\leq \mathcal{V}(\|u\|)(t+1)^{-p_{0}},~{\rm{for}}~t\geq 0,~u\in X,
\end{split}
\end{equation*}
where $\mathcal{V}$ is defined by $\mathcal{V}(r):=C\bigl(g(r)+\tilde{g}(0)\bigr),$ with $C>0$ a constant.
\end{theorem}

\par
Over the past thirty years, the coupling method has attracted considerable attention as a fundamental tool in the analysis of stochastic systems; see \cite{KS12} and the references therein for a detailed exposition. In the setting of SPDEs, coupling method has proven to be particularly effective for studying mixing properties. Notably, the coupling method is applicable to a wide range of noise types: it has been successfully employed for white noise (see \cite{KS12,M2014,Gao26JDE,N2} and references therein) as well as for bounded noise (see \cite{KNS1,Shi15,Shi21} and references therein). This versatility makes coupling method an indispensable technique in the field.

\subsection{Construct an extension}\label{Sec2.2}
First, we briefly recall some necessary preliminaries on the maximal coupling of measures, which will be used in the subsequent analysis.
We refer the reader to \cite{KS12,Shi08} for further details on this topic.

Let $X$ be a Polish space, and let $\mu,\mu'$ be two probability Borel measures on $X$. Recall that a pair $(\xi,\xi')$ of $X$-valued random variables defined on the same probability space is called a \emph{coupling} for $(\mu,\mu')$ if $\mathcal{D}(\xi)=\mu$ and $\mathcal{D}(\xi')=\mu'$.

\begin{definition}(See \cite{KS12,Shi08})
A coupling $(\xi,\xi')$ for $(\mu,\mu')$ is said to be \emph{maximal} if
\[
\mathbb{P}\{\xi\neq \xi'\}=\|\mu-\mu'\|_{\mathrm{var}}
\]
and the random variables $\xi$ and $\xi'$ conditioned on the event $\mathcal{N}:=\{\xi\neq\xi'\}$ are independent, i.e.,
\[
\mathbb{P}\{\xi\in\Gamma,\xi'\in\Gamma'\mid \mathcal{N}\}
=
\mathbb{P}\{\xi\in\Gamma\mid \mathcal{N}\}
\,
\mathbb{P}\{\xi'\in\Gamma'\mid \mathcal{N}\}
\]
for any $\Gamma,\Gamma'\in\mathcal{B}(X)$.
\end{definition}

Motivated by the ideas in \cite{KS12,N2,M2014,Shi08}, we construct an extension in the following manner. For any initial values $y,y^{\prime}\in \mathcal{H},$ let $\{y^u(t)\}_{t\geq 0}$ and $\{y^{u^{\prime}}(t)\}_{t\geq 0}$
denote the solutions of the system \eqref{WE} issued from $y,y^{\prime},$ respectively. Let $P_{N}$ denote the orthogonal projection in $H$ onto the space by the family $\{e_{1},\cdots,e_{N}\}$ with integer $N\geq 1$. With the help of $P_N$, we define an auxiliary process $\{v(t)\}_{t \ge 0}$ as the solution of the system
\begin{eqnarray}\label{F1}
 \begin{array}{l}
 \left\{
 \begin{array}{llll}
 \ddot{v}+Av+\gamma \dot{v}+f(v)+P_{N}(f(u)-f(v))=h+\eta
 \\y^v(0)=y^{\prime}
 \end{array}
 \right.
 \end{array}
 \begin{array}{lll}
 {\rm{in}}~\mathbb{R},
\\{\rm{in}}~\mathbb{R}.
\end{array}
\end{eqnarray}
Let $T>0$ be a time parameter that will be chosen later. We denote by
$\lambda_{T}(y,y^{\prime})$ and $\lambda_{T}^{\prime}(y,y^{\prime})$ the distributions of processes $\{y^v(t)\}_{0\leq t\leq T}$ and $\{y^{u^{\prime}}(t)\}_{0\leq t\leq T}$,
respectively. Then $\lambda_{T}(y,y^{\prime})$ and $\lambda_{T}^{\prime}(y,y^{\prime})$ are probability measures on $C([0,T];\mathcal{H}).$
It follows from \cite[Theorem 1.2.28]{KS12} that there exists a maximal coupling $(V_{T}(y,y^{\prime}),V_{T}^{\prime}(y,y^{\prime}))$
for the pair $(\lambda_{T}(y,y^{\prime}),\lambda_{T}^{\prime}(y,y^{\prime}))$ defined on some probability space $(\tilde{\Omega},\tilde{\mathcal{F}},\tilde{\mathbb{P}})$, the flows of this maximal coupling is defined by $\{y^{\tilde{v}}(t)\}_{0\leq t\leq T}$ and $\{y^{\tilde{u}^{\prime}}(t)\}_{0\leq t\leq T}$, respectively.
We can see that $\tilde{v}$ is the solution of the system
\begin{eqnarray}\label{29}
 \begin{array}{l}
 \left\{
 \begin{array}{llll}
 \ddot{\tilde{v}}+A\tilde{v}+\gamma\dot{\tilde{v}}+f(\tilde{v})-P_{N}f(\tilde{v})=h+\tilde{\eta}
 \\y^{\tilde{v}}(0)=y^{\prime}
 \end{array}
 \right.
 \end{array}
 \begin{array}{lll} {\rm{in}}~\mathbb{R},
\\{\rm{in}}~\mathbb{R},
\end{array}
\end{eqnarray}
where $\tilde{\eta}$ satisfies that $\mathcal{D}(\{\int_{0}^{t}\tilde{\eta}(s)ds\}_{0\leq t\leq T})=\mathcal{D}(\{\int_{0}^{t}(\eta(s)-P_{N}f(u(s))ds\}_{0\leq t\leq T})$.
Let $\tilde{u}$ be the solution of
\begin{eqnarray}\label{30}
 \begin{array}{l}
 \left\{
 \begin{array}{llll}
 \ddot{\tilde{u}}+A\tilde{u}+\gamma\dot{\tilde{u}}+f(\tilde{u})-P_{N}f(\tilde{u})=h+\tilde{\eta}
 \\y^{\tilde{u}}(0)=y
 \end{array}
 \right.
 \end{array}
 \begin{array}{lll} {\rm{in}}~\mathbb{R},
\\{\rm{in}}~\mathbb{R},
\end{array}
\end{eqnarray}
this implies that $\mathcal{D}(\{y^{\tilde{u}}(t)\}_{0\leq t\leq T})=\mathcal{D}(\{y^u(t)\}_{0\leq t\leq T})$.
We define operators $\mathcal{R}$ and $\mathcal{R}^{\prime}$ by
\begin{equation*}
\begin{split}
\mathcal{R}_{t}(y,y^{\prime},\omega)=y^{\tilde{u}}(t),~~\mathcal{R}_{t}^{\prime}(y,y^{\prime},\omega)=y^{\tilde{u}^{\prime}}(t)
\end{split}
\end{equation*}
for any $y,y^{\prime}\in \mathcal{H}, \omega\in \tilde{\Omega}, t\in [0,T]$. For each $k\ge0$, let $(\Omega^{k},\mathcal{F}^{k},\mathbb{P}^{k})$ be an independent copy of $(\tilde{\Omega},\tilde{\mathcal{F}},\tilde{\mathbb{P}})$. We then define $(\Omega,\mathcal{F},\mathbb{P})$ as the direct product of the independent probability spaces $\{(\Omega^{k},\mathcal{F}^{k},\mathbb{P}^{k})\}_{k\ge 0}$. For any $\omega=(\omega^{1},\omega^{2},\cdots)\in \Omega$
and $y,y^{\prime}\in \mathcal{H},$ we set $\tilde{y}(0)=y,\tilde{y}^{\prime}(0)=y^{\prime},$ and
\begin{equation*}
\begin{split}
&y^{\tilde{u}}(t)(\omega)=\mathcal{R}_{s}(\tilde{y}(k)(\omega),\tilde{y}^{\prime}(k)(\omega),\omega^{k}),\\
&y^{\tilde{u}^{\prime}}(t)(\omega)=\mathcal{R}_{s}^{\prime}(\tilde{y}(k)(\omega),\tilde{y}^{\prime}(k)(\omega),\omega^{k}),\\
&\mathcal{\textbf{y}}(t)=(y^{\tilde{u}}(t),y^{\tilde{u}^{\prime}}(t)),
\end{split}
\end{equation*}
where $t=s+kT, s\in[0,T).$ The above construction implies that $(\textbf{y}(t), \mathbb{P}_{\mathbf{y}})$ is an
extension for $(y^u(t), \mathbb{P}_{y}).$

\subsection{Foia\c{s}-Prodi estimate}
Foia\c{s}-Prodi estimate was firstly established in \cite{FPE}, now, it is a powerful tool to
establish the ergodicity for SPDEs, and it is often used to compensate the degeneracy of the noise, see \cite{KS12, N2} and the references therein.
The Foia\c{s}-Prodi estimate of wave equation on bounded domain was firstly established in \cite{M2014}.
Now, we establish Foia\c{s}-Prodi estimate of wave equation on the whole line.
Inspired by \cite{A1989,N2}, let us introduce a smooth space-time weight function given by
\begin{equation*}
\begin{split}
\psi(x,t):=\varphi(x)(1-\exp(-\frac{t}{\varphi(x)})), ~~(x,t)\in \mathbb{R}^{2},
\end{split}
\end{equation*}
where $\varphi(x):=\ln (x^{2}+2)$. The following properties are useful for establishing Foia\c{s}-Prodi estimate
\par
(1) $0<\psi(x,t)<\varphi(x)$ for $t>0,x\in \mathbb{R}$ and $\psi(x, 0)=0$;
\par
(2) For any $k\geq1,$ there exists some constant $C_{k}$ such that
$$|\partial^k_x\psi(x,t)|\leq C_{k}, |\partial^k_t\psi(x,t)|\leq C_{k}~{\rm{for~any}}~t\geq0, x\in \mathbb{R};$$
\par
(3) $\lim\limits_{t,|x|\rightarrow +\infty}\psi(x,t)=+\infty.$

Based on this weight function, \cite{N2} succeeded in establishing the Foia\c{s}--Prodi estimate for the white-forced complex Ginzburg--Landau equation on the whole line, and \cite{Gao26JDE} established the corresponding estimate for the white-forced hyperviscous Burgers equation on the whole line. Here, we adapt the ideas from \cite{N2,Gao26JDE,M2014} to establish the Foia\c{s}--Prodi estimate for the white-forced wave equation on the whole line, see Theorem \ref{FP}.

Let us consider the following wave equations
\begin{eqnarray}
\label{17}
&&\ddot{u}+Au+\gamma \dot{u}+f(u)=h+\eta ~~{\rm{in}}~\mathbb{R},\\
\label{18}
&&\ddot{v}+Av+\gamma \dot{v}+f(v)+P_{N}(f(u)-f(v))=h+\eta ~~{\rm{in}}~\mathbb{R}.
\end{eqnarray}
For \eqref{17} and \eqref{18}, we have the following Foia\c{s}-Prodi estimates.

\begin{theorem}\label{FP} For \eqref{17} and \eqref{18}, the following Foia\c{s}-Prodi estimates hold.
\par (1) There exist constants $\alpha>0,C>0$ such that
\begin{equation*}
\begin{split}
|y^u(t)-y^v(t)|_{\mathcal{H}}^{2}
\leq |y^u(s)-y^v(s)|_{\mathcal{H}}^{2}\exp\left (-\alpha(t-s)+C \int_{s}^{t}(\|u\|_{1}^{2}+\|v\|_{1}^{2})dr\right)
\end{split}
\end{equation*}
for any $t\geq s\geq0,$ and $N\geq 1$.
\par
(2) There exists $\alpha>0$ such that for any $\varepsilon>0,$ we can find a time $T_{0}>0$ and an integer $N\geq1$ satisfying
\begin{equation*}
\begin{split}
&|y^u(t+T_{0})-y^v(t+T_{0})|_{\mathcal{H}}^{2}\\
\leq &|y^u(s+T_{0})-y^v(s+T_{0})|_{\mathcal{H}}^{2}\exp\left (-\alpha(t-s)+C_{\ast}\varepsilon \int_{s+T_{0}}^{t+T_{0}}(\|u\|_{1}^{2}+\|v\|_{1}^{2}+\|\psi u\|_{1}^{2}+\|\psi v\|_{1}^{2})dr\right)
\end{split}
\end{equation*}
for any $t\geq s\geq0,$ where $C_{\ast}>0$ is a constant depending on $\alpha$.

\end{theorem}

\begin{proof}[Proof of Theorem \ref{FP}]
We modify the framework and arguments of \cite{N2,Gao26JDE,M2014} and adapt them to our setting. The main difference from the previous works lies in the hyperbolic nature of the wave equation. To overcome this difficulty, we employ the norm $|\cdot|_{\mathcal{H}}$.

Let $w=u-v,$ we can see that $w$ satisfies the following equation
\begin{equation}\label{FPw}
\begin{split}
\ddot{w}+Aw+\gamma \dot{w}+Q_{N}(f(u)-f(v))=0,
\end{split}
\end{equation}
where $Q_{N}:=Id-P_{N}$. It is easy to see that $y^{w}=y^{u}-y^{v}$. We multiply \eqref{FPw} by $\dot{w} + \alpha w$ and integrate by parts over the whole line $\mathbb{R}$. This yields
\begin{equation}\label{24}
\begin{split}
\partial_{t}|y^{w}|_{\mathcal{H}}^{2}=
&-2\alpha\|w\|_{1}^{2}+2(\alpha-\gamma)\|\dot{w}+\alpha w\|^{2}-2\alpha(\alpha-\gamma)(w,\dot{w}+\alpha w)\\
&-\int_{\mathbb{R}}2(\dot{w}+\alpha w)Q_{N}(f(u)-f(v))dx.
\end{split}
\end{equation}
By taking $\alpha>0$ small enough, we have
\begin{equation*}
\begin{split}
-2\alpha\|w\|_{1}^{2}+2(\alpha-\gamma)\|\dot{w}+\alpha w\|^{2}-2\alpha(\alpha-\gamma)(w,\dot{w}+\alpha w)\leq-\frac{3\alpha}{2}|y^{w}|_{\mathcal{H}}^{2}.
\end{split}
\end{equation*}
Substituting this inequality into \eqref{24}, we obtain
\begin{equation}\label{40}
\begin{split}
\partial_{t}|y^{w}|_{\mathcal{H}}^{2}\leq-\frac{3\alpha}{2}|y^{w}|_{\mathcal{H}}^{2}-\int_{\mathbb{R}}2(\dot{w}+\alpha w)Q_{N}(f(u)-f(v))dx.
\end{split}
\end{equation}
\par
(1) Noting the following fact
\begin{equation*}
\begin{split}
|2\int_{\mathbb{R}}(\dot{w}+\alpha w)Q_{N}(f(u)-f(v))dx|
\leq&2\|\dot{w}+\alpha w\|\|f(u)-f(v)\|\\
\leq& C(\|u\|_{1}^{2m}+\|v\|_{1}^{2m})\|\dot{w}+\alpha w\|\|w\|\\
\leq& C(\|u\|_{1}^{2m}+\|v\|_{1}^{2m})|y^{w}|_{\mathcal{H}}^{2}\\
\leq& \frac{1}{2}\alpha |y^{w}|_{\mathcal{H}}^{2}+ C(\alpha,m)(\|u\|_{1}^{2}+\|v\|_{1}^{2})|y^{w}|_{\mathcal{H}}^{2},
\end{split}
\end{equation*}
where the last estimate follows from applying Young's inequality to $C(\|u\|_{1}^{2m}+\|v\|_{1}^{2m})$, yielding
$C(\|u\|_{1}^{2m}+\|v\|_{1}^{2m})
\le \frac{\alpha}{2} + C(\alpha,m)\bigl(\|u\|_{1}^{2}+\|v\|_{1}^{2}\bigr)$. Substituting this inequality into \eqref{40}, we see that
\begin{equation*}
\begin{split}
\partial_{t}|y^{w}|_{\mathcal{H}}^{2}\leq-\alpha|y^{w}|_{\mathcal{H}}^{2}+C(\|u\|_{1}^{2}+\|v\|_{1}^{2})|y^{w}|_{\mathcal{H}}^{2}.
\end{split}
\end{equation*}
With the help of Gronwall inequality, we can prove (1).
\par
(2) For the proof of (2), we first recall the following lemma from \cite[Lemma 2.1]{N2}, which will play a key role in the subsequent argument.
\begin{lemma}\cite[Lemma 2.1]{N2}\label{L1}
Let us fix any $s>0.$ For any $\varepsilon,A>0,$ there is an integer $N\geq1$ such that
\begin{equation*}
\begin{split}
\|Q_{N}\chi_{A}f\|\leq \varepsilon\|f\|_{s},~~\forall f\in H^{s},
\end{split}
\end{equation*}
where $Q_{N}:=Id-P_{N}$ and $\chi_{A}$ is any smooth cut-off function from $\mathbb{R}$ to $[0,1]$
such that
\begin{equation*}
\begin{array}{l}
\chi_{A}=\left\{
\begin{array}{lll}
1\\
0
\end{array}
\right.
\end{array}
\begin{array}{lll}
x\in [-\frac{A}{2},\frac{A}{2}],
\\x\in [-A,A]^{c}.
\end{array}
\end{equation*}
\end{lemma}
In order to prove (2), we borrow some ideas from \cite{N2}. More precisely, we split $(Q_{N}(f(u)-f(v)),\dot{w}+\alpha w)$ into two parts,
\begin{equation*}
\begin{split}
&|(Q_{N}(f(u)-f(v)),\dot{w}+\alpha w)|\\
\leq&|(Q_{N}\chi_{A}(f(u)-f(v)),\dot{w}+\alpha w)|+|(Q_{N}(1-\chi_{A})(f(u)-f(v)),\dot{w}+\alpha w)|\\
=&:I_{1}+I_{2}.
\end{split}
\end{equation*}
\par
For $I_{1}$, we have
\begin{equation*}
\begin{split}
I_{1}\leq \|\dot{w}+\alpha w\|\|Q_{N}\chi_{A}(f(u)-f(v))\|
.
\end{split}
\end{equation*}
If $m\in(\frac{1}{2},1)$, we apply Lemma \ref{L1} with $s=2-2m$ to $I_{1}$, then we have
\begin{equation*}
\begin{split}
I_{1}\leq &\varepsilon\|\dot{w}+\alpha w\|\|f(u)-f(v)\|_{2-2m}
.
\end{split}
\end{equation*}
With the help of the interpolation inequality $\|\cdot\|_{2-2m}\leq C \|\cdot\|^{2m-1}\|\cdot\|_1^{2-2m}$(where $C$ is some constant depends on $m$), we have
\begin{equation*}
\begin{split}
I_{1}\leq &\varepsilon\|\dot{w}+\alpha w\|\|f(u)-f(v)\|_{2-2m}\\
\leq&\varepsilon C\|\dot{w}+\alpha w\|\|f(u)-f(v)\|^{2m-1}\|f(u)-f(v)\|_1^{2-2m}
.
\end{split}
\end{equation*}

According to the mean value theorem, by some calculations, there exists some constant $C$ which depends on $m$ such that for any $u,v\in \mathbb{R}$
\begin{equation}\label{59}
\begin{split}
|f(u)-f(v)|\leq C|u-v|(|u|^{2m}+|v|^{2m}).
\end{split}
\end{equation}
According to the property of $H^1(\mathbb{R})$ and \eqref{59}, we have
\begin{equation*}
\begin{split}
&\|f(u)-f(v)\|^{2m-1}\leq C\|w\|^{2m-1}(\|u\|_1^{2m}+\|v\|_1^{2m})^{2m-1},\\
&\|f(u)-f(v)\|_1^{2-2m}\leq C\|w\|_1^{2-2m}(\|u\|_1^{2m}+\|v\|_1^{2m})^{2-2m}
.
\end{split}
\end{equation*}
The above estimates imply that there exists some constant $C$ such that
\begin{equation*}
\begin{split}
I_{1}\leq &\varepsilon C\|\dot{w}+\alpha w\|\|w\|^{2m-1}\|w\|_{1}^{2-2m}(\|u\|_{1}^{2m}+\|v\|_{1}^{2m})\\
\leq&\varepsilon C|y^{w}|_{\mathcal{H}}^{2m}\|w\|_{1}^{2-2m}(\|u\|_{1}^{2m}+\|v\|_{1}^{2m})\\
\leq&\varepsilon C|y^{w}|_{\mathcal{H}}^{2}(\|u\|_{1}^{2m}+\|v\|_{1}^{2m})
.
\end{split}
\end{equation*}

If $m\in(0,\frac{1}{2}]$, we apply Lemma \ref{L1} with $s=1$ to $I_{1}$, then with the help of the property of $H^1(\mathbb{R})$ and \eqref{59},
we have
\begin{equation*}
\begin{split}
I_{1}&\leq \varepsilon\|\dot{w}+\alpha w\|\|f(u)-f(v)\|_{1}\\
&\leq\varepsilon C\|\dot{w}+\alpha w\|\|w\|_{1}(\|u\|_{1}^{2m}+\|v\|_{1}^{2m})\\
&\leq\varepsilon C|y^{w}|_{\mathcal{H}}\|w\|_{1}(\|u\|_{1}^{2m}+\|v\|_{1}^{2m})\\
&\leq\varepsilon C|y^{w}|_{\mathcal{H}}^{2}(\|u\|_{1}^{2m}+\|v\|_{1}^{2m})
.
\end{split}
\end{equation*}
Based on the above discussion, applying Young's inequality yields, for any $m\in(0,1)$, that
\begin{equation*}
\begin{split}
I_{1}
&\leq\varepsilon C|y^{w}|_{\mathcal{H}}^{2}(\|u\|_{1}^{2m}+\|v\|_{1}^{2m})\\
&\leq \varepsilon C|y^{w}|_{\mathcal{H}}^{2}(\|u\|^{2}_{1}+\|v\|^{2}_{1})+\frac{\alpha}{8}|y^{w}|_{\mathcal{H}}^{2}.
\end{split}
\end{equation*}
\par
For $I_{2}$, it follows from the property of the weight function $\psi$ that
\begin{equation*}
\begin{split}
\psi^{-2m}(x,t)\leq \varepsilon,~~\forall t\geq T_{0}, |x|\geq\frac{A}{2},
\end{split}
\end{equation*}
provided that $T_{0},A>0$ are large enough. This implies that
\begin{equation*}
\begin{split}
I_{2}&=((1-\chi_{A})(f(u)-f(v)),Q_{N}(\dot{w}+\alpha w))\\
&=(Q_{N}(\dot{w}+\alpha w),(1-\chi_{A})(|\psi u|^{2m}u-|\psi v|^{2m}v)\psi^{-2m})\\
&\leq \varepsilon\|\dot{w}+\alpha w\|\||\psi u|^{2m}u-|\psi v|^{2m}v\|.
\end{split}
\end{equation*}
It follows from the property of $H^1(\mathbb{R})$ and \eqref{59} that
\begin{equation*}
\begin{split}
\||\psi u|^{2m}u-|\psi v|^{2m}v\|
&=\||\psi|^{2m}(f(u)-f(v))\|\\
&\leq C\||\psi|^{2m}|u-v|(|u|^{2m}+|v|^{2m})\|\\
&= C\||w|(|\psi u|^{2m}+|\psi v|^{2m})\|\\
&\leq  C\|w\|(\|\psi u\|_{1}^{2m}+\|\psi v\|_{1}^{2m}),
\end{split}
\end{equation*}
this implies that
\begin{equation*}
\begin{split}
I_{2}
&\leq \varepsilon C\|\dot{w}+\alpha w\|\|w\|(\|\psi u\|_{1}^{2m}+\|\psi v\|_{1}^{2m})\\
&\leq \varepsilon C|y^{w}|_{\mathcal{H}}^{2}(\|\psi u\|^{2}_{1}+\|\psi v\|^{2}_{1})+\frac{\alpha}{8}|y^{w}|_{\mathcal{H}}^{2},
\end{split}
\end{equation*}
where the last estimate follows from applying Young's inequality to $C(\|\psi u\|_{1}^{2m}+\|\psi v\|_{1}^{2m})$, yielding
$C(\|\psi u\|_{1}^{2m}+\|\psi v\|_{1}^{2m})
\le \frac{\alpha}{8} + C(\alpha,m)\bigl(\|\psi u\|_{1}^{2}+\|\psi v\|_{1}^{2}\bigr)$. This leads to
\begin{equation*}
\begin{split}
\partial_{t}|y^{w}|_{\mathcal{H}}^{2}+[\alpha-\varepsilon C(\|u\|_{1}^{2}+\|v\|_{1}^{2}+\|\psi u\|_{1}^{2}+\|\psi v\|_{1}^{2})]|y^{w}|_{\mathcal{H}}^{2}\leq0.
\end{split}
\end{equation*}
With the help of Gronwall inequality, we can prove (2).
\par
This completes the proof.
\end{proof}

\section{Weighted energy estimates for the stochastic wave equation}
\subsection{Energy identity for linear stochastic wave equation}
Let us first consider the linear stochastic wave equation
\begin{eqnarray}\label{LWE}
 \begin{array}{l}
 \left\{
 \begin{array}{llll}
\ddot{u}+Au+\gamma \dot{u}=f+\eta
 \\u(x,0)=u_{0}(x)
 \\\dot{u}(x,0)=u_{1}(x)
 \end{array}
 \right.
 \end{array}
 \begin{array}{lll}
 {\rm{in}}~\mathbb{R}\times(0,+\infty),\\
 {\rm{in}}~\mathbb{R},\\
 {\rm{in}}~\mathbb{R},
\end{array}
\end{eqnarray}
where $f=f(t)$ is a continuous predictable $H-$valued process and $\eta$ is given exactly as in \eqref{37}.
\par
We introduce the following functionals
\begin{align*}
&|y^{u}|_{\mathcal{H}}^{2}:=\|u\|^{2}+\|u_{x}\|^{2}+\|\dot{u}+\alpha u\|^{2},\\
&|y^{u}|_{\psi}^{2}:=\|\psi u\|^{2}+\|\psi u_{x}\|^{2}+\|\psi (\dot{u}+\alpha u)\|^{2},\\
&\mathcal{E}_{u}(t):=|y^{u}|_{\mathcal{H}}^{2}+\frac{1}{m+1}\int_{\mathbb{R}}|u|^{2m+2}dx,\\
&\mathcal{E}_{u}^{\psi}(t):=\mathcal{E}_{u}(t)+|y^{u}|_{\psi}^{2}+\frac{1}{m+1}\int_{\mathbb{R}}\psi^{2}|u|^{2m+2}dx,\\
&\mathcal{F}_{u}(t):=\mathcal{E}_{u}(t)+\alpha\int_{0}^{t}\mathcal{E}_{u}(s)ds,\\
&\mathcal{F}_{u}^{\psi}(t):=\mathcal{E}_{u}^{\psi}(t)+\alpha\int_{0}^{t}\mathcal{E}_{u}^{\psi}(s)ds,\\
&\mathcal{F}_{u}^{p}(t):=(\mathcal{E}_{u}(t))^{p}+\alpha p\int_{0}^{t}(\mathcal{E}_{u}(s))^{p}ds,\\
&\mathcal{F}_{u}^{\psi,p}(t):=(\mathcal{E}_{u}^{\psi}(t))^{p}+\alpha p\int_{0}^{t}(\mathcal{E}_{u}^{\psi}(s))^{p}ds,
\end{align*}
where $\alpha>0$ will be chosen sufficiently small in the later analysis.
\begin{proposition}\label{Keypro}
For any solution $u$ of \eqref{LWE}, we have
\begin{equation}\label{E1}
\begin{split}
d|y^{u}|_{\mathcal{H}}^{2}=&[-2\alpha\|u\|_{1}^{2}+2(\alpha-\gamma)\|\dot{u}+\alpha u\|^{2}-2\alpha(\alpha-\gamma)(u,\dot{u}+\alpha u)\\
&+\int_{\mathbb{R}}2(\dot{u}+\alpha u)fdx+\mathcal{B}_{1}]dt+d\bar{M}(t)
,
\end{split}
\end{equation}
where $d\bar{M}(t):=2(\dot{u}+\alpha u,dW(t))$. Moreover, the following weight equality holds
\begin{equation}\label{E2}
\begin{split}
d|y^{u}|_{\psi}^{2}
=&[-2\alpha\int_{\mathbb{R}}\psi^{2}(u^{2}+u_{x}^{2})dx+2(\alpha-\gamma)\|\psi (\dot{u}+\alpha u)\|^{2}-2\alpha(\alpha-\gamma)(\psi u,\psi (\dot{u}+\alpha u))\\
&+\int_{\mathbb{R}}2\psi\psi_{t}u^{2}dx+\int_{\mathbb{R}}2\psi\psi_{t}(\dot{u}+\alpha u)^{2}dx+\int_{\mathbb{R}}2\psi\psi_{t}u_{x}^{2}dx
-\int_{\mathbb{R}}4\psi\psi_{x}(\dot{u}+\alpha u)u_{x}dx\\
&+\int_{\mathbb{R}}2\psi^{2} (\dot{u}+\alpha u)fdx+\sum_{j=1}^{\infty}b_{j}^{2}\|\psi e_{j}\|^{2}]dt+d\bar{M}^{\psi}(t)
,
\end{split}
\end{equation}
where $d\bar{M}^{\psi}(t):=2(\psi^{2}(\dot{u}+\alpha u),dW(t)).$ In particular, for any solution $u$ of \eqref{LWE} with $\eta\equiv0$ (the deterministic linear wave equation), it holds that
\begin{equation}\label{E3}
\begin{split}
\frac{d}{dt}|y^{u}|_{\mathcal{H}}^{2}=-2\alpha\|u\|_{1}^{2}+2(\alpha-\gamma)\|\dot{u}+\alpha u\|^{2}-2\alpha(\alpha-\gamma)(u,\dot{u}+\alpha u)+2((\dot{u}+\alpha u),f)
.
\end{split}
\end{equation}
\end{proposition}
Its proof can be found in Appendix.
\begin{remark}
Proposition~\ref{Keypro} serves as a fundamental tool in the derivation of the estimates for solutions to system \eqref{WE}.
\end{remark}

\subsection{Probability estimates for solutions}
\begin{proposition}\label{pro5}
Let the assumption $\textbf{(A)}$ hold. Then, there exist some constants $\alpha>0,C>0$ such that, for any $t\geq0, T\geq0,$ the following a priori estimate holds for solutions of \eqref{WE}
\begin{equation}\label{2}
\begin{split}
\mathcal{E}_{u}(t+T)-\mathcal{E}_{u}(T)+\alpha\int_{T}^{t+T}\mathcal{E}_{u}(s)ds
\leq (M(t)-\frac{\beta}{2}\langle M(t)\rangle)+C(\|h\|^{2}+\mathcal{B}_{1})t,
\end{split}
\end{equation}
where $\beta:=\frac{\alpha}{8\mathcal{B}_{1}}$ and $M(t)$ is defined by
$$M(t):=2\int_{T}^{t+T}(\dot{u}+\alpha u,dW(s))=2\sum\limits_{j=1}^{\infty}\int_{T}^{t+T}b_{j}(\dot{u}+\alpha u,e_{j})d\beta_{j}(s).$$
Moreover, we have the probability estimates
\begin{eqnarray}
\label{3}&&\mathbb{P}\{\sup\limits_{t\geq 0}\left(\mathcal{E}_{u}(t+T)+\int_{T}^{t+T}\alpha\mathcal{E}_{u}(s)ds-C(\|h\|^{2}+\mathcal{B}_{1})t\right)\geq\mathcal{E}_{u}(T)+\rho\}
\leq e^{-\beta \rho},\\
\label{12}&&\mathbb{P}\{\sup\limits_{t\geq 0}\left(\mathcal{F}_{u}(t+T)-C(\|h\|^{2}+\mathcal{B}_{1})t\right)\geq\mathcal{F}_{u}(T)+\rho\}
\leq e^{-\beta \rho},
\end{eqnarray}
for all $\rho>0$.
\end{proposition}
\begin{proof}[Proof of Proposition \ref{pro5}]
By taking $\alpha>0$ small enough, we have
\begin{equation}\label{25}
\begin{split}
-2\alpha\int_{\mathbb{R}}(u^{2}+u_{x}^{2})dx+2(\alpha-\gamma)\|\dot{u}+\alpha u\|^{2}-2\alpha(\alpha-\gamma)(u,\dot{u}+\alpha u)\leq
-\frac{7}{4}\alpha|y|_{\mathcal{H}}^{2}.
\end{split}
\end{equation}
If we take $f=-|u|^{2m}u+h$ in \eqref{E1}, then
\begin{equation*}
\begin{split}
2\int_{\mathbb{R}}(\dot{u}+\alpha u)fdx=-\frac{1}{m+1}\frac{d}{dt}\int_{\mathbb{R}}|u|^{2m+2}dx-2\alpha\int_{\mathbb{R}}|u|^{2m+2}dx+2\int_{\mathbb{R}}(\dot{u}+\alpha u)hdx
,
\end{split}
\end{equation*}
this implies that
\begin{equation}\label{49}
\begin{split}
d\mathcal{E}_{u}(t)
=&[-2\alpha\|u\|_{1}^{2}+2(\alpha-\gamma)\|\dot{u}+\alpha u\|^{2}-2\alpha(\alpha-\gamma)(u,\dot{u}+\alpha u)
-2\alpha\int_{\mathbb{R}}|u|^{2m+2}dx]dt
\\&+2(\dot{u}+\alpha u,h)dt+\mathcal{B}_{1}dt+d\bar{M}(t)
.
\end{split}
\end{equation}
\eqref{25} gives that
\begin{equation}\label{19}
\begin{split}
d\mathcal{E}_{u}(t)
&\leq(-\frac{3}{2}\alpha|y|_{\mathcal{H}}^{2}
-2\alpha\int_{\mathbb{R}}|u|^{2m+2}dx)dt+d\bar{M}(t)+C(\|h\|^{2}+\mathcal{B}_{1})dt\\
&\leq-\frac{3}{2}\alpha(|y|_{\mathcal{H}}^{2}+\frac{1}{m+1}\int_{\mathbb{R}}|u|^{2m+2}dx)dt+d\bar{M}(t)+C(\|h\|^{2}+\mathcal{B}_{1})dt
.
\end{split}
\end{equation}
Integrating \eqref{19} over $[T, t+T]$, we obtain
\begin{equation*}
\begin{split}
\mathcal{E}_{u}(t+T)-\mathcal{E}_{u}(T)
&\leq\int_{T}^{t+T}(-\frac{3}{2}\alpha|y|_{\mathcal{H}}^{2}-2\alpha\int_{\mathbb{R}}|u|^{2m+2}dx)ds
+M(t)+C(\|h\|^{2}+\mathcal{B}_{1})t\\
&\leq-\frac{3}{2}\alpha\int_{T}^{t+T}(|y|_{\mathcal{H}}^{2}+\frac{1}{m+1}\int_{\mathbb{R}}|u|^{2m+2}dx)ds
+M(t)+C(\|h\|^{2}+\mathcal{B}_{1})t\\
&=-\frac{3}{2}\alpha\int_{T}^{t+T}\mathcal{E}_{u}(s)ds
+M(t)-\frac{\beta}{2}\langle M(t)\rangle+\frac{\beta}{2}\langle M(t)\rangle+C(\|h\|^{2}+\mathcal{B}_{1})t,
\end{split}
\end{equation*}
where
\begin{equation*}
\begin{split}
\langle M(t)\rangle
=4\sum_{j=1}^{\infty} b_{j}^{2}\int_{T}^{t+T}(\dot{u}+\alpha u,e_{j})^{2}ds\leq 4\mathcal{B}_{1}\sum_{j=1}^{\infty}\int_{T}^{t+T}(\dot{u}+\alpha u,e_{j})^{2}ds
=4\mathcal{B}_{1}\int_{T}^{t+T}\|\dot{u}+\alpha u\|^{2}ds.
\end{split}
\end{equation*}
By taking $\beta=\frac{\alpha}{8\mathcal{B}_{1}}$, we have
$
\beta\langle M(t)\rangle\leq \frac{\alpha}{2} \int_{T}^{t+T}|y|_{\mathcal{H}}^{2}ds,
$
thus, we have \eqref{2}.
It follows from \eqref{2} that
\begin{equation*}
\begin{split}
\mathcal{E}_{u}(t+T)+\int_{T}^{t+T}[\alpha\mathcal{E}_{u}(s)-C(\|h\|^{2}+\mathcal{B}_{1})]ds-\mathcal{E}_{u}(T)
\leq M(t)-\frac{\beta}{2}\langle M(t)\rangle,
\end{split}
\end{equation*}
this implies that
\begin{equation*}
\begin{split}
\{\sup\limits_{t\geq 0}\left(\mathcal{E}_{u}(t+T)+\int_{T}^{t+T}[\alpha\mathcal{E}_{u}(s)-(\|h\|^{2}+\mathcal{B}_{1})]ds-\mathcal{E}_{u}(T)\right)\geq \rho\}
\subset
\{\sup\limits_{t\geq 0}\left(M(t)-\frac{\beta}{2}\langle M(t)\rangle\right)\geq \rho\},
\end{split}
\end{equation*}
with the help of the exponential supermartingale inequality, we have \eqref{3}.
\par
Noting the fact
\begin{equation*}
\begin{split}
\mathcal{F}_{u}(t+T)-\mathcal{F}_{u}(T)
=\mathcal{E}_{u}(t+T)-\mathcal{E}_{u}(T)+\int_{T}^{t+T}\alpha\mathcal{E}_{u}(s)ds,
\end{split}
\end{equation*}
together with \eqref{3}, we can prove \eqref{12}.

This completes the proof.
\end{proof}

\begin{proposition}\label{pro3}
Let the assumption $\textbf{(A)}$ hold. Then, there exists a constant $\alpha>0$ such that, for any $p\geq2, q>1$,
there exist positive constants $C, L_{p}$ such that
\begin{equation*}
\begin{split}
\mathbb{P}(\sup\limits_{t\geq 0}[\mathcal{F}^{p}_{u}(t+T)-\mathcal{F}^{p}_{u}(T)-(L_{p}+1)t-2]\geq \rho)\leq CQ_{q}(\rho+1),
\end{split}
\end{equation*}
for all $\rho>0,T\geq 0,$ and $y_{0}\in \mathcal{H}$, where $C$ is a constant depending on $y_{0},p,q,h,\mathcal{B}_{1}$.
\end{proposition}
\begin{proof}[Proof of Proposition \ref{pro3}]
By applying Ito formula to $\mathcal{E}_{u}^{p}(t)$, we have
\begin{equation*}
\begin{split}
d\mathcal{E}_{u}^{p}(t)
=p\mathcal{E}_{u}^{p-1}(t)d\mathcal{E}_{u}(t)+\frac{1}{2}p(p-1)\mathcal{E}_{u}^{p-2}(t)(d\mathcal{E}_{u}(t))^{2}
,
\end{split}
\end{equation*}
this gives us
\begin{equation*}
\begin{split}
d\mathcal{E}_{u}^{p}(t)
=&p\mathcal{E}_{u}^{p-1}(t)[(-2\alpha\|u\|_{1}^{2}+2(\alpha-\gamma)\|\dot{u}+\alpha u\|^{2}-2\alpha(\alpha-\gamma)(u,\dot{u}+\alpha u)
-2\alpha\int_{\mathbb{R}}|u|^{2m+2}dx)dt
\\&~~~~~~~~~~~~~~+2(\dot{u}+\alpha u,h)dt+d\bar{M}(t)+\mathcal{B}_{1}dt]+\frac{1}{2}p(p-1)\mathcal{E}_{u}^{p-2}(t)(d\mathcal{E}_{u}(t))^{2}
.
\end{split}
\end{equation*}
Since $(d\mathcal{E}_{u}(t))^{2}=\|\dot{u}+\alpha u\|^{2}\mathcal{B}_{1}dt$ and with the help of \eqref{25}, we have
\begin{equation*}
\begin{split}
&d\mathcal{E}_{u}^{p}(t)\\
=&p\mathcal{E}_{u}^{p-1}(t)[(-2\alpha\|u\|_{1}^{2}+2(\alpha-\gamma)\|\dot{u}+\alpha u\|^{2}-2\alpha(\alpha-\gamma)(u,\dot{u}+\alpha u)
-2\alpha\int_{\mathbb{R}}|u|^{2m+2}dx)dt
\\&~~~~~~~~~~~~~~+2(\dot{u}+\alpha u,h)dt+d\bar{M}(t)+\mathcal{B}_{1}dt]+\frac{1}{2}p(p-1)\mathcal{E}_{u}^{p-2}(t)\|\dot{u}+\alpha u\|^{2}\mathcal{B}_{1}dt
\\\leq & p\mathcal{E}_{u}^{p-1}(t)[(-\frac{3}{2}\alpha|y|_{\mathcal{H}}^{2}-2\alpha\int_{\mathbb{R}}|u|^{2m+2}dx+C\|h\|^{2}+\mathcal{B}_{1})dt+d\bar{M}(t)]
+\frac{1}{2}p(p-1)\mathcal{E}_{u}^{p-2}(t)\|\dot{u}+\alpha u\|^{2}\mathcal{B}_{1}dt
\\\leq & -\frac{3}{2}\alpha p\mathcal{E}_{u}^{p}(t)dt+C(p,h,\mathcal{B}_{1})\mathcal{E}_{u}^{p-1}(t)dt+C\mathcal{E}_{u}^{p-1}(t)d\bar{M}(t)
\\\leq & -\alpha p\mathcal{E}_{u}^{p}(t)dt+C(p,h,\mathcal{B}_{1},\alpha)dt+C\mathcal{E}_{u}^{p-1}(t)d\bar{M}(t)
.
\end{split}
\end{equation*}
Namely, it holds that
\begin{equation}\label{47}
\begin{split}
d\mathcal{E}_{u}^{p}(t)
\leq-\alpha p\mathcal{E}_{u}^{p}(t)dt+C(p,h,\mathcal{B}_{1},\alpha)dt+C\mathcal{E}_{u}^{p-1}(t)d\bar{M}(t)
.
\end{split}
\end{equation}
By applying Gronwall inequality, we have
\begin{equation}\label{55}
\begin{split}
\mathbb{E}\mathcal{E}_{u}^{p}(t)\leq e^{-\alpha pt}\mathbb{E}\mathcal{E}_{u}^{p}(0)+C(p,h,\mathcal{B}_{1},\alpha)~~\forall t\geq0.
\end{split}
\end{equation}
It follows from \eqref{47} that there exists some positive constant $L_{p}$ such that
\begin{equation*}
\begin{split}
\mathcal{F}_{u}^{p}(t+T)-\mathcal{F}_{u}^{p}(T)-(L_{p}+1)t-2
\leq M^{p}(t)-t-2
,
\end{split}
\end{equation*}
where $M^{p}(t):=C\int_{T}^{t+T}\mathcal{E}_{u}^{p-1}(s)d\bar{M}(s)$.
\par
We define $M^{p}_{*}(t):=\sup_{s\in [0,t]}|M^{p}(s)|.$
For any $q>1$, according to \eqref{55}, the Burkholder-Davis-Gundy and H\"{o}lder inequalities, we have
\begin{equation}\label{48}
\begin{split}
\mathbb{E}(M^{p}_{*}(t))^{2q}
&\leq C\mathbb{E}\langle M^{p}(t)\rangle^{q}
\leq C\mathbb{E}(\int_{T}^{t+T}\mathcal{E}_{u}^{2p-1}(s)ds)^{q}\\
&\leq Ct^{q-1}\mathbb{E}\int_{T}^{t+T}\mathcal{E}_{u}^{(2p-1)q}ds
\leq C(u_{0},u_{1},p,q,h,\mathcal{B}_{1})(t+1)^{q}.
\end{split}
\end{equation}
Noting the following fact for all $\rho>0,$ it holds that
\begin{equation*}
\begin{split}
\{\sup\limits_{t\geq 0}[M^{p}(t)-t-2]\geq \rho\}
\subset &\bigcup\limits_{m\geq 0}\{\sup\limits_{t\in [m,m+1)}[M^{p}(t)-t-2]\geq \rho\}\\
\subset &\bigcup\limits_{m\geq 0}\{M_{*}^{p}(m+1)\geq \rho+m+2\}.
\end{split}
\end{equation*}
For any $q>1$, it follows from the above analysis, the maximal martingale inequality, the Chebyshev inequality and \eqref{48} that
\begin{equation*}
\begin{split}
\mathbb{P}\{\sup\limits_{t\geq 0}[M^{p}(t)-t-2]\geq \rho\}
&\leq\mathbb{P}(\bigcup\limits_{m\geq 0}\{M_{*}^{p}(m+1)\geq \rho+m+2\})\\
&\leq\sum\limits_{m\geq 0}\mathbb{P}(M_{*}^{p}(m+1)\geq \rho+m+2)\\
&\leq\sum\limits_{m\geq 0}\frac{\mathbb{E}(M_{*}^{p}(m+1))^{2q}}{(\rho+m+2)^{2q}}\\
&\leq C(y_{0},p,q,h,\mathcal{B}_{1})\sum\limits_{m\geq 0}\frac{(m+2)^{q}}{(\rho+m+2)^{2q}}\\
&\leq C(y_{0},p,q,h,\mathcal{B}_{1})\sum\limits_{m\geq 0}\frac{1}{(\rho+m+2)^{q}}\\
&\leq \frac{C(y_{0},p,q,h,\mathcal{B}_{1})}{(\rho+1)^{q-1}}.
\end{split}
\end{equation*}
\par
This completes the proof.
\end{proof}

\begin{proposition}\label{pro2}
Let the assumption $\textbf{(A)}$ hold. Then, there exist some constants $\alpha>0,C>0$ such that, for any $t\geq0, T\geq0,$ the following a priori estimate holds for solutions of \eqref{WE}
\begin{equation}\label{4}
\begin{split}
&\mathcal{E}_{u}^{\psi}(t+T)+\alpha\int_{T}^{t+T}\mathcal{E}_{u}^{\psi}(s)ds
-C(\mathcal{B}_{1}+\|h\|^{2}+\mathcal{B}_{2}+\|\varphi h\|^{2})t-\mathcal{E}_{u}^{\psi}(T)-C\mathcal{E}_{u}(T)\\
\leq &C(M(t)-\frac{\beta_{0}}{2}\langle M(t)\rangle)+(M^{\psi}(t)-\frac{\beta_{0}}{2}\langle M^{\psi}(t)\rangle)
,
\end{split}
\end{equation}
where $\beta_{0}:=\frac{\alpha}{8\mathcal{B}_{1}}\wedge\frac{\alpha}{8\mathcal{B}_{2}}$ and $M^{\psi}(t)$ is defined by
$$M^{\psi}(t):=2\int_{T}^{t+T}(\psi^{2}(\dot{u}+\alpha u),dW(s))=2\sum\limits_{j=1}^{\infty}\int_{T}^{t+T}b_{j}(\psi(\dot{u}+\alpha u),\psi e_{j})d\beta_{j}(s).$$
In other words, we have
\begin{equation}\label{16}
\begin{split}
&\mathcal{F}_{u}^{\psi}(t+T)-\mathcal{F}_{u}^{\psi}(T)-C(\mathcal{B}_{1}+\|h\|^{2}+\mathcal{B}_{2}+\|\varphi h\|^{2})t-C\mathcal{E}_{u}(T)\\
\leq &C(M(t)-\frac{\beta_{0}}{2}\langle M(t)\rangle)+(M^{\psi}(t)-\frac{\beta_{0}}{2}\langle M^{\psi}(t)\rangle)
.
\end{split}
\end{equation}
Moreover, there exist some positive constants $\mathcal{K},\mathcal{M},\hat{\beta}$ such that
\begin{eqnarray}
\label{6}
&&\mathbb{P}\{\sup\limits_{t\geq 0}\left(\mathcal{E}^{\psi}_{u}(t+T)+\int_{T}^{t+T}\alpha\mathcal{E}^{\psi}_{u}(s)ds-\mathcal{K}t\right)\geq\mathcal{E}_{u}^{\psi}(T)
+\mathcal{M}\mathcal{E}_{u}(T)+\rho\}\leq 2e^{-\hat{\beta} \rho},\\
\label{7}
&&\mathbb{P}\{\sup\limits_{t\geq 0}\left(\mathcal{F}^{\psi}_{u}(t+T)-\mathcal{K}t\right)\geq\mathcal{F}^{\psi}_{u}(T)+\mathcal{M}\mathcal{E}_{u}(T)+\rho\}
\leq 2e^{-\hat{\beta}\rho},
\end{eqnarray}
for all $\rho>0$.
\end{proposition}
\begin{proof}[Proof of Proposition \ref{pro2}]
By taking $\alpha>0$ small enough, we have
\begin{equation*}
\begin{split}
&-2\alpha\int_{\mathbb{R}}\psi^{2}(u^{2}+u_{x}^{2})dx+2(\alpha-\gamma)\|\psi (\dot{u}+\alpha u)\|^{2}-2\alpha(\alpha-\gamma)(\psi u,\psi (\dot{u}+\alpha u))
\leq-\frac{7}{4}\alpha|y|_{\psi}^{2},\\
&\int_{\mathbb{R}}[2\psi\psi_{t}u^{2}+2\psi\psi_{t}(\dot{u}+\alpha u)^{2}+2\psi\psi_{t}u_{x}^{2}-4\psi\psi_{x}(\dot{u}+\alpha u)u_{x}]dx
\leq\frac{\alpha}{8}|y|_{\psi}^{2}+C|y|_{\mathcal{H}}^{2}
\leq\frac{\alpha}{8}|y|_{\psi}^{2}+C\mathcal{E}_{u}(t).
\end{split}
\end{equation*}
If we take $f=-|u|^{2m}u+h$ in \eqref{E2}, then
\begin{equation*}
\begin{split}
&2\int_{\mathbb{R}}\psi^{2}(\dot{u}+\alpha u)fdx\\
=&-\frac{1}{m+1}\frac{d}{dt}\int_{\mathbb{R}}\psi^{2}|u|^{2m+2}dx
+\int_{\mathbb{R}}[\frac{2}{m+1}\psi\psi_{t}|u|^{2m+2}-2\alpha\psi^{2}|u|^{2m+2}+2\psi^{2}(\dot{u}+\alpha u)h]dx.
\end{split}
\end{equation*}
This gives us
\begin{equation}\label{46}
\begin{split}
&d(|y^{u}|_{\psi}^{2}+\frac{1}{m+1}\int_{\mathbb{R}}\psi^{2}|u|^{2m+2}dx)\\
=&[-2\alpha\int_{\mathbb{R}}\psi^{2}(u^{2}+u_{x}^{2})dx+2(\alpha-\gamma)\|\psi (\dot{u}+\alpha u)\|^{2}-2\alpha(\alpha-\gamma)(\psi u,\psi (\dot{u}+\alpha u)-2\alpha\int_{\mathbb{R}}\psi^{2}|u|^{2m+2}dx\\
&+\int_{\mathbb{R}}[\frac{2}{m+1}\psi\psi_{t}|u|^{2m+2}+2\psi\psi_{t}u^{2}+2\psi\psi_{t}(\dot{u}+\alpha u)^{2}+2\psi\psi_{t}u_{x}^{2}
-4\psi\psi_{x}(\dot{u}+\alpha u)u_{x}]dx\\
&+\int_{\mathbb{R}}2\psi^{2} (\dot{u}+\alpha u)h dx+\sum_{j=1}^{\infty}b_{j}^{2}\|\psi e_{j}\|^{2}]dt+d\bar{M}^{\psi}(t)
.
\end{split}
\end{equation}
Noting the fact
\begin{equation*}
\begin{split}
\frac{2}{m+1}\int_{\mathbb{R}}\psi\psi_{t}|u|^{2m+2}dx
&\leq\frac{\alpha}{100(m+1)}\int_{\mathbb{R}}\psi^{2}|u|^{2m+2}dx+C\int_{\mathbb{R}}|u|^{2m+2}dx\\
&\leq\frac{\alpha}{100(m+1)}\int_{\mathbb{R}}\psi^{2}|u|^{2m+2}dx+C\mathcal{E}_{u}(t)
,
\end{split}
\end{equation*}
we have
\begin{equation}\label{20}
\begin{split}
&d(|y^{u}|_{\psi}^{2}+\frac{1}{m+1}\int_{\mathbb{R}}\psi^{2}|u|^{2m+2}dx)\\
\leq&-\frac{3\alpha}{2}(|y|_{\psi}^{2}+\frac{1}{m+1}\int_{\mathbb{R}}\psi^{2}|u|^{2m+2}dx)dt
+2\sum\limits_{j=1}^{\infty}b_{j}(\psi^{2} (\dot{u}+\alpha u),e_{j})d\beta_{j}\\
&+C\mathcal{E}_{u}(t)dt+C(\sum\limits_{i=1}^{\infty}b_{i}^{2}\|\psi e_{i}\|^{2}+\|\psi h\|^{2})dt.
\end{split}
\end{equation}
With the help of \eqref{19} and \eqref{20}, we have
\begin{equation}\label{21}
\begin{split}
&\mathcal{E}_{u}^{\psi}(t+T)-\mathcal{E}_{u}^{\psi}(T)+\frac{3}{2}\alpha\int_{T}^{t+T}\mathcal{E}_{u}^{\psi}(s)ds\\
\leq &M(t)+ M^{\psi}(t)+C\int_{T}^{t+T}\mathcal{E}_{u}(s)ds+\int_{T}^{t+T}(\mathcal{B}_{1}+\|h\|^{2}+\sum\limits_{i=1}^{\infty}b_{i}^{2}\|\psi e_{i}\|^{2}+\|\psi h\|^{2})ds\\
\leq &M(t)+ M^{\psi}(t)+C\int_{T}^{t+T}\mathcal{E}_{u}(s)ds+C(\mathcal{B}_{1}+\|h\|^{2}+\mathcal{B}_{2}+\|\varphi h\|^{2})t.
\end{split}
\end{equation}
Noting the fact
\begin{equation*}
\begin{split}
\langle M^{\psi}\rangle(t)=4\sum\limits_{j=1}^{\infty}\int_{T}^{t+T}b_{j}^{2}(\psi(\dot{u}+\alpha u),\psi e_{j})^{2}ds
\leq4\mathcal{B}_{2}\int_{T}^{t+T}\|\psi(\dot{u}+\alpha u)\|^{2}ds
\end{split}
\end{equation*}
and by taking $\beta_{0}=\frac{\alpha}{8\mathcal{B}_{1}}\wedge\frac{\alpha}{8\mathcal{B}_{2}}$, we have
\begin{equation*}
\begin{split}
\beta_{0}\langle M(t)\rangle\leq \frac{\alpha}{2} \int_{T}^{t+T}|y|_{\mathcal{H}}^{2}ds,~~
\beta_{0}\langle M^{\psi}(t)\rangle\leq \frac{\alpha}{2} \int_{T}^{t+T}|y|_{\psi}^{2}ds.
\end{split}
\end{equation*}
This implies that
\begin{equation}\label{5}
\begin{split}
&\mathcal{E}_{u}^{\psi}(t+T)-\mathcal{E}_{u}^{\psi}(T)+\alpha\int_{T}^{t+T}\mathcal{E}_{u}^{\psi}(s)ds-C(\mathcal{B}_{1}+\|h\|^{2}+\mathcal{B}_{2}+\|\varphi h\|^{2})t\\
\leq &(M(t)-\frac{\beta_{0}}{2}\langle M(t)\rangle)+(M^{\psi}(t)-\frac{\beta_{0}}{2}\langle M^{\psi}(t)\rangle)+C\int_{T}^{t+T}\mathcal{E}_{u}(s)ds.
\end{split}
\end{equation}
\par
It follows from \eqref{5} and \eqref{2} that for some suitable constants $C,C_1,C_2,C_3$
\begin{equation*}
\begin{split}
&\mathcal{E}_{u}^{\psi}(t+T)+\alpha\int_{T}^{t+T}\mathcal{E}_{u}^{\psi}(s)ds
-C(\mathcal{B}_{1}+\|h\|^{2}+\mathcal{B}_{2}+\|\varphi h\|^{2})t-\mathcal{E}_{u}^{\psi}(T)-C\mathcal{E}_{u}(T)\\
\leq &(M(t)-\frac{\beta_{0}}{2}\langle M(t)\rangle)+(M^{\psi}(t)-\frac{\beta_{0}}{2}\langle M^{\psi}(t)\rangle)+C_1(\alpha\int_{T}^{t+T}\mathcal{E}_{u}(s)ds-\mathcal{E}_{u}(T)-C_2(\mathcal{B}_{1}+\|h\|^{2})t)\\
\leq &C_3(M(t)-\frac{\beta_{0}}{2}\langle M(t)\rangle)+(M^{\psi}(t)-\frac{\beta_{0}}{2}\langle M^{\psi}(t)\rangle),
\end{split}
\end{equation*}
this implies that \eqref{4} holds. According to \eqref{4} and the definition of $\mathcal{F}_{u}^{\psi}$,
we can obtain \eqref{16}. With the help of the exponential supermartingale inequality, \eqref{4} and \eqref{16}, we have \eqref{6}
and \eqref{7}.

This completes the proof.
\end{proof}
\begin{proposition}\label{pro4}
Let the assumption $\textbf{(A)}$ hold. Then, there exists a constant $\alpha>0$ such that, for any $p\geq2, q>1$, there exist positive constants $C, L_{p}$ such that
\begin{equation*}
\begin{split}
\mathbb{P}(\sup\limits_{t\geq 0}[\mathcal{F}^{\psi,p}_{u}(t+T)-\mathcal{F}^{\psi,p}_{u}(T)-(L_{p}+1)t-2]\geq \rho)\leq CQ_{q}(\rho+1),
\end{split}
\end{equation*}
for all $\rho>0,T\geq 0,$ and $y_{0}\in \mathcal{H}$, where $C$ is a constant depending on $y_{0},p,q,h,\mathcal{B}_{1},\mathcal{B}_{2}$.
\end{proposition}
\begin{proof}[Proof of Proposition \ref{pro4}]
It follows from \eqref{49} and \eqref{46} that
\begin{align*}
&d\mathcal{E}_{u}^{\psi}(t)\\
=&[-2\alpha\|u\|_{1}^{2}+2(\alpha-\gamma)\|\dot{u}+\alpha u\|^{2}-2\alpha(\alpha-\gamma)(u,\dot{u}+\alpha u)
-2\alpha\int_{\mathbb{R}}|u|^{2m+2}dx\\
&-2\alpha\int_{\mathbb{R}}\psi^{2}(u^{2}+u_{x}^{2})dx+2(\alpha-\gamma)\|\psi (\dot{u}+\alpha u)\|^{2}-2\alpha(\alpha-\gamma)(\psi u,\psi (\dot{u}+\alpha u)-2\alpha\int_{\mathbb{R}}\psi^{2}|u|^{2m+2}dx\\
&+\int_{\mathbb{R}}[\frac{2}{m+1}\psi\psi_{t}|u|^{2m+2}+2\psi\psi_{t}u^{2}+2\psi\psi_{t}(\dot{u}+\alpha u)^{2}+2\psi\psi_{t}u_{x}^{2}
-4\psi\psi_{x}(\dot{u}+\alpha u)u_{x}]dx\\
&+\int_{\mathbb{R}}2\psi^{2} (\dot{u}+\alpha u)hdx+2(\dot{u}+\alpha u,h)+\sum_{j=1}^{\infty}b_{j}^{2}\|\psi e_{j}\|^{2}+\mathcal{B}_{1}]dt+d\bar{M}(t)+d\bar{M}^{\psi}(t)\\
=&:P(t)dt+d\bar{\mathcal{M}}(t),
\end{align*}
where $d\bar{\mathcal{M}}(t):=d\bar{M}(t)+d\bar{M}^{\psi}(t)=2((1+\psi^{2})(\dot{u}+\alpha u),dW(t))$.
By applying Ito formula to $(\mathcal{E}_{u}^{\psi}(t))^{p}$, we have
\begin{equation}\label{51}
\begin{split}
d(\mathcal{E}_{u}^{\psi}(t))^{p}
&=p(\mathcal{E}_{u}^{\psi}(t))^{p-1}d\mathcal{E}_{u}^{\psi}(t)+\frac{1}{2}p(p-1)(\mathcal{E}_{u}^{\psi}(t))^{p-2}(d\mathcal{E}_{u}^{\psi}(t))^{2}\\
&=p(\mathcal{E}_{u}^{\psi}(t))^{p-1}P(t)dt+p(\mathcal{E}_{u}^{\psi}(t))^{p-1}d\bar{\mathcal{M}}(t)+\frac{1}{2}p(p-1)(\mathcal{E}_{u}^{\psi}(t))^{p-2}d\langle \bar{\mathcal{M}}\rangle(t)
.
\end{split}
\end{equation}
This implies that for $t\geq0,$ it holds that
\begin{equation*}
\begin{split}
&\mathbb{E}(\mathcal{E}_{u}^{\psi}(t))^{p}\\
=&\mathbb{E}(\mathcal{E}_{u}^{\psi}(0))^{p}+p\int_{0}^{t}\mathbb{E}(\mathcal{E}_{u}^{\psi}(s))^{p-1}P(s)ds
+\frac{1}{2}p(p-1)\int_{0}^{t}\mathbb{E}(\mathcal{E}_{u}^{\psi}(s))^{p-2}d\langle \bar{\mathcal{M}}\rangle(s)\\
=&\mathbb{E}(\mathcal{E}_{u}^{\psi}(0))^{p}+p\int_{0}^{t}\mathbb{E}(\mathcal{E}_{u}^{\psi}(s))^{p-1}P(s)ds
+2p(p-1)\int_{0}^{t}\mathbb{E}(\mathcal{E}_{u}^{\psi}(s))^{p-2}\sum\limits_{j=1}^{\infty}b_{j}^{2}((1+\psi^{2})(\dot{u}+\alpha u),e_{j})^{2}ds.
\end{split}
\end{equation*}
where we have used $d\langle \bar{\mathcal{M}}\rangle(s)=4\sum\limits_{j=1}^{\infty}b_{j}^{2}((1+\psi^{2})(\dot{u}+\alpha u),e_{j})^{2}ds$.
This yields that
\begin{equation}\label{50}
\begin{split}
\frac{d}{dt}\mathbb{E}(\mathcal{E}_{u}^{\psi}(t))^{p}
&=p\mathbb{E}(\mathcal{E}_{u}^{\psi}(t))^{p-1}P(t)
+2p(p-1)\mathbb{E}(\mathcal{E}_{u}^{\psi}(t))^{p-2}\sum\limits_{j=1}^{\infty}b_{j}^{2}((1+\psi^{2})(\dot{u}+\alpha u),e_{j})^{2}.
\end{split}
\end{equation}

By taking $\alpha>0$ small enough and Young inequality, we have the following facts
\begin{equation*}
\begin{split}
&P(s)\leq -\frac{3}{2}\alpha\mathcal{E}_{u}^{\psi}(s)+C(p,h,\psi,\mathcal{B}_{1},\mathcal{B}_{2}),\\
&\sum\limits_{j=1}^{\infty}b_{j}^{2}((1+\psi^{2})(\dot{u}+\alpha u),e_{j})^{2}
\leq C(\mathcal{B}_{1},\mathcal{B}_{2})(\|\dot{u}+\alpha u\|^{2}+\|\psi(\dot{u}+\alpha u)\|^{2})\leq C(\mathcal{B}_{1},\mathcal{B}_{2})\mathcal{E}_{u}^{\psi}(t)
.
\end{split}
\end{equation*}
Putting the above facts together with \eqref{50} and using Gronwall's inequality gives
\begin{equation*}
\begin{split}
\mathbb{E}(\mathcal{E}_{u}^{\psi}(t))^{p}\leq e^{-\alpha pt}\mathbb{E}(\mathcal{E}_{u}^{\psi}(0))^{p}+C(\alpha,p,h,\psi,\mathcal{B}_{1},\mathcal{B}_{2})~~\forall t\geq0.
\end{split}
\end{equation*}
It follows from \eqref{51} that there exists some positive constant $L_{p}$ such that
\begin{equation*}
\begin{split}
\mathcal{F}_{u}^{\psi,p}(t+T)-\mathcal{F}_{u}^{\psi,p}(T)-(L_{p}+1)t-2
\leq \mathcal{M}_{p}(t)-t-2
,
\end{split}
\end{equation*}
where $\mathcal{M}_{p}(t):=C\int_{T}^{t+T}(\mathcal{E}_{u}^{\psi}(s))^{p-1}d\bar{\mathcal{M}}(s)$.
\par
The rest of the proof is similar to that of Proposition \ref{pro3}, we omit it.
\end{proof}

For any $p\geq1$, we define the stopping times
\begin{equation*}
\begin{split}
\tau^{u}_{p}:=\inf\{t\geq0:\mathcal{F}_{u}^{\psi,p}(t)\geq M\mathcal{E}_{u}^{p}(0)+(K+L)t+\rho\},
\end{split}
\end{equation*}
where $M,K,L,\rho\geq 0$.

\begin{proposition}\label{pro6}
Let the assumption $\textbf{(A)}$ hold. Then, there exists a constant $C>0$ such that if $K\geq \mathcal{K}, M\geq 1+\mathcal{M},$
we have
\begin{equation}\label{56}
\begin{split}
\mathbb{P}(l\leq\tau^{u}_{1}<+\infty)\leq Ce^{-\beta_{0}(\rho+Ll) },
\end{split}
\end{equation}
for all $L,l\geq0,$ and $y_{0}\in \mathcal{H},$ where $C$ is a constant depending on $y_{0},h,\mathcal{B}_{1},\mathcal{B}_{2}$.
\end{proposition}
\begin{proof}[Proof of Proposition \ref{pro6}]
According to the definition of $\psi$, we have $\psi(x,0)=0$. It follows from the definitions of $\mathcal{E}_{u}$ and $\mathcal{F}_{u}^{\psi}$ that
$\mathcal{F}_{u}^{\psi}(0)=\mathcal{E}_{u}^{\psi}(0)= \mathcal{E}_{u}(0).$ On the event $\{l\leq\tau^{u}_{1}<\infty\},$ the definition of $\tau^{u}_{1}$ implies that
\begin{equation*}
\begin{split}
\mathcal{F}_{u}^{\psi}(\tau^{u}_{1})&\geq M\mathcal{E}_{u}(0)+(K+L)\tau^{u}_{1}+\rho\\
&\geq M\mathcal{E}_{u}(0)+\mathcal{K}\tau^{u}_{1}+Ll+\rho\\
&\geq \mathcal{F}_{u}^{\psi}(0)+\mathcal{M}\mathcal{E}_{u}(0)+\mathcal{K}\tau^{u}_{1}+Ll+\rho,
\end{split}
\end{equation*}
thus, we have
\begin{equation*}
\begin{split}
\sup\limits_{t\geq 0}\left(\mathcal{F}^{\psi}_{u}(t)-\mathcal{F}_{u}^{\psi}(0)-\mathcal{M}\mathcal{E}_{u}(0)-\mathcal{K}t\right)\geq Ll+\rho,
\end{split}
\end{equation*}
with the help of \eqref{7}, we have \eqref{56}.
\end{proof}

\begin{proposition}\label{pro1}
Let the assumption $\textbf{(A)}$ hold. Then, for any $p\geq2,q>1$, there exists a positive constant $C$ such that if $K\geq \mathcal{K}, M\geq 1+\mathcal{M},$
we have
\begin{equation}\label{14}
\begin{split}
\mathbb{P}(l\leq\tau^{u}_{p}<+\infty)\leq CQ_{q}(\rho+Ll+1),
\end{split}
\end{equation}
for all $L,l\geq0,$ and $y_{0}\in \mathcal{H},$ where $C$ is a constant depending on $y_{0},p,q,h,\mathcal{B}_{1},\mathcal{B}_{2}$.
\end{proposition}
\begin{proof}[Proof of Proposition \ref{pro1}]
According to the definition of $\psi$, we have $\psi(x,0)=0$. It follows from the definitions of $\mathcal{E}_{u}^{p}$ and $\mathcal{F}_{u}^{\psi,p}$ that
$\mathcal{F}_{u}^{\psi,p}(0)=(\mathcal{E}_{u}^{\psi}(0))^{p}= \mathcal{E}_{u}^{p}(0).$
On the event $\{l\leq\tau^{u}_{p}<\infty\},$ the definition of $\tau^{u}_{p}$ implies that
\begin{equation*}
\begin{split}
\mathcal{F}_{u}^{\psi,p}(\tau^{u}_{p})&\geq M\mathcal{E}_{u}^{p}(0)+(K+L)\tau^{u}_{p}+\rho\\
&\geq M\mathcal{E}_{u}^{p}(0)+\mathcal{K}\tau^{u}_{p}+Ll+\rho\\
&\geq \mathcal{F}_{u}^{\psi,p}(0)+\mathcal{M}\mathcal{E}_{u}^{p}(0)+\mathcal{K}\tau^{u}_{p}+Ll+\rho,
\end{split}
\end{equation*}
thus, we have
\begin{equation*}
\begin{split}
\sup\limits_{t\geq 0}\left(\mathcal{F}^{\psi,p}_{u}(t)-\mathcal{F}_{u}^{\psi,p}(0)-\mathcal{M}\mathcal{E}_{u}^{p}(0)-\mathcal{K}t\right)\geq Ll+\rho,
\end{split}
\end{equation*}
with the help of Proposition \ref{pro4}, we have \eqref{14}.
\end{proof}

We define the stopping time $\tau^{u}$ by
$$\tau^{u}:=\tau^{u}_{1}\wedge \tau^{u}_{2}.$$
\begin{proposition}\label{pro0}
Let the assumption $\textbf{(A)}$ hold. Then, for any $q>1$, there exists a positive constant $C$ such that if $K\geq \mathcal{K}, M\geq 1+\mathcal{M},$
we have
\begin{equation}\label{57}
\begin{split}
\mathbb{P}(l\leq\tau^{u}<+\infty)\leq CQ_{q}(\rho+Ll+1),
\end{split}
\end{equation}
for all $L,l\geq0,$ and $y_{0}\in \mathcal{H},$ where $C$ is a constant depending on $y_{0},q,h,\mathcal{B}_{1},\mathcal{B}_{2}$.
\end{proposition}
\begin{proof}[Proof of Proposition \ref{pro0}]
We only need to note the fact
$$\{l\leq\tau^{u}<+\infty\}\subset\{l\leq\tau^{u}_{1}<+\infty\}\cup \{l\leq\tau^{u}_{2}<+\infty\}.$$
With the help of Proposition \ref{pro6} and Proposition \ref{pro1}, we can prove \eqref{57}.
\end{proof}

\subsection{Probability estimates for the truncated process}
Let $\{u(t)\}_{t\ge 0}$ and $\{u'(t)\}_{t\ge 0}$ denote the solutions of system \eqref{WE} with initial conditions $y\in\mathcal{H}$ and $y'\in\mathcal{H}$, respectively, and let $\{v(t)\}_{t\ge 0}$ denote the solution of system \eqref{F1} starting from $y'\in\mathcal{H}$.
Let $$\tau=\tau^{v}\wedge\tau^{u}\wedge\tau^{u^{\prime}}.$$ We introduce the truncated processes $\{\hat{u}(t)\}_{t\geq0},$ $\{\hat{u}^{\prime}(t)\}_{t\geq0}$
and $\{\hat{v}(t)\}_{t\geq0}$ as follows:
\begin{equation*}
\begin{split}
\hat{u}(t)=
\left\{
\begin{array}{lll}
u(t),~&{\rm{for}}~0\leq t\leq \tau,\\
{\rm{solves}}~ \ddot{z}+Az+\gamma \dot{z}+|z|^{2m}z=0,~&{\rm{for}}~t\geq \tau,
\end{array}
\right.
\\
\hat{u}^{\prime}(t)=
\left\{
\begin{array}{lll}
u^{\prime}(t),~&{\rm{for}}~0\leq t\leq \tau,\\
{\rm{solves}}~ \ddot{z}+Az+\gamma \dot{z}+|z|^{2m}z=0,~&{\rm{for}}~t\geq \tau,
\end{array}
\right.
\\
\hat{v}(t)=
\left\{
\begin{array}{lll}
v(t),~&{\rm{for}}~0\leq t\leq \tau,\\
{\rm{solves}}~ \ddot{z}+Az+\gamma \dot{z}+|z|^{2m}z=0,~&{\rm{for}}~t\geq \tau.
\end{array}
\right.
\end{split}
\end{equation*}

\begin{proposition}\label{proM2}
Let the assumption $\textbf{(A)}$ hold. Then, for any $p\geq2,q>1$, there exist positive constants $C, \hat{K}$ and $\hat{M}$ such that
\begin{equation*}
\begin{split}
\mathbb{P}\{\sup\limits_{t\geq 0}\left(\mathcal{F}^{\psi,p}_{\hat{u}^{\prime}}(t)-\hat{K}t\right)\geq \hat{M}\mathcal{E}_{u^{\prime}}^{p}(0)+\rho\}
\leq CQ_{q}(\rho+1),
\end{split}
\end{equation*}
for all $\rho>0,$ and $y,y^{\prime}\in \mathcal{H},$ where $C$ is a constant depending on $y,y^{\prime},p,q,h,\mathcal{B}_{1},\mathcal{B}_{2}$.
\end{proposition}
\begin{proof}[Proof of Proposition \ref{proM2}]
Let us first consider the auxiliary equation
\begin{equation*}
\begin{split}
\ddot{z}+Az+\gamma \dot{z}+|z|^{2m}z=0~~{\rm{in}}~\mathbb{R}.
\end{split}
\end{equation*}
With the help of the similar arguments as in \eqref{47} and \eqref{51}, we have
\begin{equation*}
\begin{split}
d\mathcal{F}_{z}^{p}(t)\leq 0,~~d\mathcal{F}_{z}^{\psi,p}(t)\leq C\mathcal{E}_{z}^{p}(t)dt.
\end{split}
\end{equation*}
Namely, for any $t\geq 0,T\geq0$, we have
\begin{equation*}
\begin{split}
\mathcal{F}_{z}^{p}(t+T)\leq \mathcal{F}_{z}^{p}(T),~~\mathcal{F}_{z}^{\psi,p}(t+T)\leq \mathcal{F}_{z}^{\psi,p}(T)+C\int_{T}^{t+T}\mathcal{E}_{z}^{p}(s)ds.
\end{split}
\end{equation*}
The above results imply that
$
\mathcal{F}_{z}^{\psi,p}(t+T)\leq C\mathcal{F}_{z}^{\psi,p}(T).
$
This implies that on the event $\{\tau<\infty\},$ it holds that
\begin{equation*}
\begin{split}
\mathcal{F}_{\hat{u}^{\prime}}^{\psi,p}(t+\tau)\leq C\mathcal{F}_{u^{\prime}}^{\psi,p}(\tau)~~{\rm{for}}~t\geq0.
\end{split}
\end{equation*}
\par
We have the fact there exists a constant $C^{\prime}>1$ such that
\begin{equation}\label{22}
\begin{split}
\mathcal{F}_{\hat{u}^{\prime}}^{\psi,p}(t)-C^{\prime}(L_{p}+1)t\leq C^{\prime}\sup\limits_{t\geq0}(\mathcal{F}_{u^{\prime}}^{\psi,p}(t)-(L_{p}+1)t)~~{\rm{for}}~t\geq0.
\end{split}
\end{equation}
Indeed, \eqref{22} is obvious on the event $\{\tau<\infty\}$. On the event $\{\tau=+\infty\}$, it is easy to see that
\begin{equation*}
\begin{split}
\mathcal{F}_{\hat{u}^{\prime}}^{\psi,p}(t)-C^{\prime}(L_{p}+1)t\leq \sup\limits_{t\geq0}(\mathcal{F}_{u^{\prime}}^{\psi,p}(t)-(L_{p}+1)t)~~{\rm{for}}~t\geq0.
\end{split}
\end{equation*}
Thus, \eqref{22} implies that for $t\geq0$,
\begin{equation*}
\begin{split}
C^{\prime}\sup\limits_{t\geq0}(\mathcal{F}_{u^{\prime}}^{\psi,p}(t)-(L_{p}+1)t-\mathcal{E}_{u^{\prime}}^{p}(0))
\geq\mathcal{F}_{\hat{u}^{\prime}}^{\psi,p}(t)-C^{\prime}(L_{p}+1)t-C^{\prime}\mathcal{E}_{u^{\prime}}^{p}(0)
.
\end{split}
\end{equation*}
By taking $\hat{K}=C^{\prime}(L_{p}+1), \hat{M}=C^{\prime},$ with the help of Proposition \ref{pro4}, we can prove Proposition \ref{proM2}.
\par
This completes the proof.
\end{proof}

\begin{proposition}\label{proT2}
Let the assumption $\textbf{(A)}$ hold. Then, for any $p\geq2,q>1$, there exists a positive constant $C$ such that if $K\geq \hat{K}$ and $M\geq \hat{M}$
\begin{equation}\label{15}
\begin{split}
\mathbb{P}(l\leq\tau^{\hat{u}^{\prime}}_{p}<+\infty)\leq CQ_{q}(\rho+Ll+1),
\end{split}
\end{equation}
for all $L,l\geq0,\rho>0,$ and $y,y^{\prime}\in \mathcal{H},$ where $C$ is a constant depending on $y,y^{\prime},p,q,h,\mathcal{B}_{1},\mathcal{B}_{2}$.
\end{proposition}
\begin{proof}[Proof of Proposition \ref{proT2}]
On the event $\{l\leq\tau^{\hat{u}^{\prime}}_{p}<\infty\},$ the definition of $\tau^{\hat{u}^{\prime}}_{p}$ implies that
\begin{equation*}
\begin{split}
\mathcal{F}^{\psi,p}_{\hat{u}^{\prime}}(\tau^{\hat{u}^{\prime}}_{p})
\geq M\mathcal{E}_{u^{\prime}}^{p}(0)+(K+L)\tau^{\hat{u}^{\prime}}_{p}+\rho
\geq \hat{M}\mathcal{E}_{u^{\prime}}^{p}(0)+\hat{K}\tau^{\hat{u}^{\prime}}_{p}+Ll+\rho,
\end{split}
\end{equation*}
thus, we have
\begin{equation*}
\begin{split}
\sup\limits_{t\geq 0}\left(\mathcal{F}^{\psi,p}_{\hat{u}^{\prime}}(\tau^{\hat{u}^{\prime}}_{p})-\hat{K}t-\hat{M}\mathcal{E}_{u^{\prime}}^{p}(0)\right)\geq Ll+\rho,
\end{split}
\end{equation*}
with the help of Proposition \ref{proM2}, we have \eqref{15}.
\end{proof}
By repeating the arguments in Propositions~\ref{pro6} and~\ref{proM2}, and with the help of Proposition~\ref{pro2}, we obtain the following result.
\begin{proposition}\label{pro7}
Let the assumption $\textbf{(A)}$ hold. Then, there exist positive constants $C$ and $\hat{\beta}$ such that if $K\geq \mathcal{K}, M\geq 1+\mathcal{M},$
we have
\begin{equation}
\begin{split}
\mathbb{P}(l\leq\tau^{\hat{u}^{\prime}}_{1}<+\infty)\leq Ce^{-\hat{\beta}(\rho+Ll)},
\end{split}
\end{equation}
for all $L,l\geq0,\rho>0,$ and $y,y^{\prime}\in \mathcal{H},$ where $C$ is a constant depends on $y,y^{\prime},h,\mathcal{B}_{1},\mathcal{B}_{2}$.
\end{proposition}
\par
According to Proposition \ref{proT2} and Proposition \ref{pro7}, we have the following result.
\begin{proposition}\label{pro01}
Let the assumption $\textbf{(A)}$ hold. Then, for any $q>1$, there exists a positive constant $C$ such that if $K\geq \mathcal{K}, M\geq 1+\mathcal{M},$
we have
\begin{equation}
\begin{split}
\mathbb{P}(l\leq\tau^{\hat{u}^{\prime}}<+\infty)\leq CQ_{q}(\rho+Ll+1),
\end{split}
\end{equation}
for all $L,l\geq0,$ and $y,y^{\prime}\in \mathcal{H},$ where $C$ is a constant depending on $y,y^{\prime},q,h,\mathcal{B}_{1},\mathcal{B}_{2}$.
\end{proposition}

\subsection{Probability estimate for the auxiliary process}
To handle the integral terms appearing in the Foia\c{s}-Prodi estimate,
we establish an estimate for the stopping time $\tau^{v}$.
\begin{proposition}\label{pro10}
Let the assumption $\textbf{(A)}$ hold. Then, for any $q>1$, if \eqref{41} holds for some $N\geq 1$, then there exists a positive constant $C$ such that
\begin{equation}\label{60}
\begin{split}
\mathbb{P}(\tau^{v}<\infty)
\leq CQ_{q}(\rho+1)+\frac{1}{2}\left(\exp(Ce^{C(\rho+\mathcal{E}_{u}^{2}(0)+\mathcal{E}_{v}^{2}(0))}d^{2})-1 \right)^{\frac{1}{2}}
\end{split}
\end{equation}
for any $\rho>0$ and $y,y^{\prime}\in \mathcal{H}$ with $d=|y-y^{\prime}|_{\mathcal{H}}.$
\end{proposition}

\begin{proof}[Proof of Proposition \ref{pro10}]
In order to prove Proposition~\ref{pro10}, we borrow some ideas from \cite{M2014,KS12,Shi08,N2,Gao26JDE} and adapt them to our setting.
Let us define vectors $\hat{e}_{j}=(0,e_j)$ and their vector span
$$\mathcal{H}_{N}:=span\{\hat{e}_{1},\hat{e}_{2},\cdots,\hat{e}_{N}\},$$
which is an $N-$dimensional subspace of $\mathcal{H}.$
Let $\mathcal{P}_{N}$ be the orthogonal projection from $\mathcal{H}$ to $\mathcal{H}_{N}$.
Without loss of generality, we take the underlying probability space $(\Omega,\mathcal{F},\mathbb{P})$ to be of the following form: let $\Omega:=C_{0}([0,+\infty);\mathcal{H})$ be the space of continuous $\mathcal{H}$-valued paths vanishing at $t=0$, let $\mathbb{P}$ be the law of the Wiener process
$$
\hat{\xi}(t)=\sum_{i=1}^{\infty} b_i \beta_i(t)\hat{e}_i,
$$
and let $\mathcal{F}$ be the completion of the Borel $\sigma$-algebra with respect to the topology of uniform convergence on compact sets.

The proof is divided into the following steps.
\par
\textbf{Step 1}. We now derive an estimate for $\mathbb{P}(\tau^{v}<\infty)$.

Applying Proposition \ref{pro0} with $l=0$ yields
\begin{equation}\label{62}
\begin{split}
\mathbb{P}(\tau^{v}<\infty)\leq CQ_{q}(\rho+1)+\mathbb{P}(\tau^{\hat{v}}<\infty).
\end{split}
\end{equation}

We now provide an estimate for $\mathbb{P}(\tau^{\hat{v}}<\infty)$. For any integer $N\geq1,$ we define the transform
\begin{equation}\label{58}
\begin{split}
\Phi^{u,u^{\prime}}:&~\Omega\longrightarrow \Omega\\
&\omega(t)\rightarrow \omega(t)-\int_{0}^{t}\mathbb{I}_{s\leq \tau}\mathcal{P}_{N}(0,f(\hat{u})-f(\hat{v}))ds.
\end{split}
\end{equation}
Thanks to the pathwise uniqueness for the stochastic wave equation, it follows that
\begin{equation*}
\begin{split}
\mathbb{P}\{y^{\hat{u}^{\prime}}(\Phi^{u,u^{\prime}}(\omega),t)=y^{\hat{v}}(\omega,t),~\forall t\geq0\}=1,
\end{split}
\end{equation*}
this leads to
\begin{equation}\label{61}
\begin{split}
\mathbb{P}(\tau^{\hat{v}}<\infty)=\Phi^{u,u^{\prime}}_{\ast}\mathbb{P}(\tau^{\hat{u}^{\prime}}<\infty)\leq \mathbb{P}(\tau^{\hat{u}^{\prime}}<\infty)+\|\mathbb{P}-\Phi^{u,u^{\prime}}_{\ast}\mathbb{P}\|_{var}.
\end{split}
\end{equation}
Substituting \eqref{61} into \eqref{62} yields
\begin{equation}\label{63}
\begin{split}
\mathbb{P}(\tau^{v}<\infty)\leq CQ_{q}(\rho+1)+\mathbb{P}(\tau^{\hat{u}^{\prime}}<\infty)+\|\mathbb{P}-\Phi^{u,u^{\prime}}_{\ast}\mathbb{P}\|_{var}.
\end{split}
\end{equation}
\par
\textbf{Step 2}. To estimate the total variation distance $\|\mathbb{P}-\Phi^{u,u'}_{\ast}\mathbb{P}\|_{\mathrm{var}}$, we begin by estimating $y^{\hat{u}}(t)-y^{\hat{v}}(t)$. For $0\le t\le \tau$, it holds that
\begin{equation}\label{11}
\begin{split}
|y^{\hat{u}}(t)-y^{\hat{v}}(t)|_{\mathcal{H}}^{2}\leq C\exp\left (-\frac{a}{2}t+C(\rho+\mathcal{E}_{u}^{2}(0)+\mathcal{E}_{v}^{2}(0))\right)d^{2}.
\end{split}
\end{equation}
\par
Indeed, noting the fact $y^{\hat{u}}(t)-y^{\hat{v}}(t)=y^{u}(t)-y^{v}(t)$ for $0\leq t\leq\tau$ and the definition of $\tau$,
we have
\begin{equation}\label{9}
\begin{split}
&\mathcal{F}_{\hat{u}}^{\psi,2}(t)=\mathcal{F}_{u}^{\psi,2}(t)\leq M\mathcal{E}_{u}^{2}(0)+(K+L)t+\rho,\\
&\mathcal{F}_{\hat{v}}^{\psi,2}(t)=\mathcal{F}_{v}^{\psi,2}(t)\leq M\mathcal{E}_{v}^{2}(0)+(K+L)t+\rho,
\end{split}
\end{equation}
for $0\leq t\leq\tau$. Let $T_{0}$ be the number in (2) of Theorem \ref{FP}. Now, we distinguish the following two cases.
\par
\textit{Case 1. $\tau\leq T_{0}.$}
\par
Foia\c{s}-Prodi estimate in (1) of Theorem \ref{FP} and \eqref{9} show that
\begin{equation*}
\begin{split}
|y^{u}(t)-y^{v}(t)|_{\mathcal{H}}^{2}&\leq |y^{u}(0)-y^{v}(0)|_{\mathcal{H}}^{2}\exp\left (-\alpha t+C\int_{0}^{t}(\|u\|_{1}^{2}+\|v\|_{1}^{2})dr\right)\\
&\leq C|y^{u}(0)-y^{v}(0)|_{\mathcal{H}}^{2}\exp\left (-\alpha t+C[(K+L)T_{0}+\mathcal{E}_{u}^{2}(0)+\mathcal{E}_{v}^{2}(0)+\rho]\right),
\end{split}
\end{equation*}
this implies that
\begin{equation}\label{8}
\begin{split}
|y^{u}(t)-y^{v}(t)|_{\mathcal{H}}^{2}\leq C\exp\left(-\alpha t+C(\rho+\mathcal{E}_{u}^{2}(0)+\mathcal{E}_{v}^{2}(0))\right)d^{2}.
\end{split}
\end{equation}
\par
\textit{Case 2. $\tau> T_{0}.$}
\par
We can see that \eqref{8} holds on $[0,T_{0}]$. We apply Foia\c{s}-Prodi estimate in (2) of  Theorem \ref{FP} with $\varepsilon=\frac{\alpha}{4C_{\ast}(K+L)}$ to $y^{u}(t)-y^{v}(t)$ on $[T_{0},\tau)$, it holds that
\begin{equation*}
\begin{split}
|y^{u}(t)-y^{v}(t)|_{\mathcal{H}}^{2}\leq |y^{u}(T_{0})-y^{v}(T_{0})|_{\mathcal{H}}^{2}\exp\left (-\alpha(t-T_{0})+C_{\ast}\varepsilon \int_{T_{0}}^{t}(\|u\|_{1}^{2}+\|v\|_{1}^{2}+\|\psi u\|_{1}^{2}+\|\psi v\|_{1}^{2})dr\right).
\end{split}
\end{equation*}
\eqref{9} leads to
\begin{equation*}
\begin{split}
|y^{u}(t)-y^{v}(t)|_{\mathcal{H}}^{2}\leq C\exp\left (-\frac{\alpha}{2}t+C(\rho+\mathcal{E}_{u}^{2}(0)+\mathcal{E}_{v}^{2}(0))\right)d^{2}.
\end{split}
\end{equation*}
Case 1 and Case 2 imply \eqref{11}.
\par
\textbf{Step 3}. The estimate of $\|\mathbb{P}-\Phi^{u,u^{\prime}}_{\ast}\mathbb{P}\|_{var}$.
\par
The space $\Omega:=C_{0}([0,+\infty);\mathcal{H})$ can be represented in the form
$$\Omega=\Omega_{N}\oplus\Omega_{N}^{\bot},$$
where $\Omega_{N}=C([0,+\infty);\mathcal{H}_{N})$ and $\Omega_{N}^{\bot}=C([0,+\infty);\mathcal{H}_{N}^{\bot})$. Let $\mathbb{P}_{N}$ and $\mathbb{P}_{N}^{\bot}$ be the images of
$\mathbb{P}$ under the natural projections $\Pi_{N}:\Omega\rightarrow\Omega_{N}$ and $\Pi^{N}:\Omega\rightarrow\Omega_{N}^{\bot}$, respectively.
For $\omega=(\omega^{(1)},\omega^{(2)})\in \Omega,$
\begin{equation*}
\begin{split}
\Phi^{u,u^{\prime}}(\omega)=\Phi^{u,u^{\prime}}(\omega^{(1)},\omega^{(2)})=(\Psi^{u,u^{\prime}}(\omega^{(1)},\omega^{(2)}),\omega^{(2)}),
\end{split}
\end{equation*}
where
\begin{equation*}
\begin{split}
\Psi^{u,u^{\prime}}:&~\Omega\longrightarrow C([0,+\infty);\mathcal{H}_{N})\\
&\omega(t)\rightarrow \omega^{(1)}(t)+\int_{0}^{t}\mathcal{A}(s,\omega^{(1)},\omega^{(2)})ds
\end{split}
\end{equation*}
with $\mathcal{A}(s):=-\mathbb{I}_{s\leq \tau}\mathcal{P}_{N}(0,f(\hat{u})-f(\hat{v})).$ By applying \cite[Lemma 3.3.13]{KS12}, we have
\begin{equation*}
\begin{split}
\|\mathbb{P}-\Phi^{u,u^{\prime}}_{\ast}\mathbb{P}\|_{var}\leq \int_{C([0,+\infty);Q_{N}H)}\|\mathbb{P}_{N}-\Psi^{u,u^{\prime}}_{\ast}(\mathbb{P}_{N},\omega^{(2)})\|_{var}\mathbb{P}_{N}^{\bot}(d\omega^{(2)}).
\end{split}
\end{equation*}
Now, we prove the following Novikov condition: If for any $C>0$ and $\omega^{(2)},$ it holds that
\begin{equation}\label{10}
\begin{split}
\mathbb{E}_{N}\exp(C\int_{0}^{\infty}\|\mathcal{A}(t)\|^{2}dt)<\infty,
\end{split}
\end{equation}
\par
Indeed, since $\|\mathcal{A}(t)\|\leq C\cdot \mathbb{I}_{t\leq \tau}\cdot\|w\|_{1}(\|u\|_{1}^{2}+\|v\|_{1}^{2}+1)$ and \eqref{11}, we have
\begin{equation*}
\begin{split}
\int_{0}^{\infty}\|\mathcal{A}(t)\|^{2}dt
&=\int_{0}^{\tau}\|\mathcal{A}(t)\|^{2}dt\\
&\leq C\int_{0}^{\tau}\|w\|^{2}_{1}(\|u\|_{1}^{4}+\|v\|_{1}^{4}+1)dt\\
&\leq C\int_{0}^{\tau}|y^{u}(t)-y^{v}(t)|_{\mathcal{H}}^{2}(\|u\|_{1}^{4}+\|v\|_{1}^{4}+1)dt\\
&\leq C|y^{u}(0)-y^{v}(0)|_{\mathcal{H}}^{2}\int_{0}^{\tau}\exp\left (-\frac{\alpha}{2}t+C(\rho+\mathcal{E}_{u}^{2}(0)+\mathcal{E}_{v}^{2}(0))\right)(\|u\|_{1}^{4}+\|v\|_{1}^{4}+1)dt
.
\end{split}
\end{equation*}
It follows from \eqref{9} that
\begin{equation*}
\begin{split}
\int_{0}^{\infty}\|\mathcal{A}(t)\|^{2}dt
\leq C|y^{u}(0)-y^{v}(0)|_{\mathcal{H}}^{2}\int_{0}^{\tau}\exp\left (-\frac{\alpha}{2}t\right) H(t)dt
,
\end{split}
\end{equation*}
where $H(t):=\exp\left (C(\rho+\mathcal{E}_{u}^{2}(0)+\mathcal{E}_{v}^{2}(0))\right)[\rho+M\mathcal{E}_{u}^{2}(0)+M\mathcal{E}_{v}^{2}(0)+(K+L)t+1],$
then, we have
\begin{equation*}
\begin{split}
\int_{0}^{\infty}\|\mathcal{A}(t)\|^{2}dt
\leq Ce^{C(\rho+\mathcal{E}_{u}^{2}(0)+\mathcal{E}_{v}^{2}(0))}d^{2}
,
\end{split}
\end{equation*}
namely, we prove Novikov condition \eqref{10}. With the help of Novikov condition \eqref{10}, Girsanov theorem shows that
\begin{equation*}
\begin{split}
\|\mathbb{P}_{N}-\Psi^{u,u^{\prime}}_{\ast}(\mathbb{P}_{N},\omega^{(2)})\|_{var}
\leq \frac{1}{2}[(\mathbb{E}_{N}\exp(6\sup_{1\leq i\leq N}b_{i}^{-2}\int_{0}^{\infty}\|\mathcal{A}(t)\|^{2}dt))^{\frac{1}{2}}-1]^{\frac{1}{2}},
\end{split}
\end{equation*}
here we recall that, by assumption \eqref{41}, $b_i\neq 0$ for $i=1,2,\dots,N$.
Then, we can obtain
\begin{equation}\label{13}
\begin{split}
\|\mathbb{P}-\Phi^{u,u^{\prime}}_{\ast}\mathbb{P}\|_{var}\leq \frac{1}{2}\left(\exp(Ce^{C(\rho+\mathcal{E}_{u}^{2}(0)+\mathcal{E}_{v}^{2}(0))}d^{2})-1 \right)^{\frac{1}{2}}.
\end{split}
\end{equation}
\par
\textbf{Step 4}. Proof of estimate \eqref{60}.
\par
Combining \eqref{63}, \eqref{13}, and the estimate for $\mathbb{P}(\tau^{\hat{u}'}<\infty)$ from Proposition~\ref{proT2} yields \eqref{60}.

This completes the proof.
\end{proof}

\section{Irreducibility for the stochastic wave equation}

In this section, we establish the irreducibility of the stochastic wave equation \eqref{WE} within the framework developed in \cite{KS12,Gao22ECOCV}. We first prove the continuity of the solution operator associated with the deterministic wave equation; more precisely, we establish the continuity of the mapping from noise to solution along controllers (see Proposition \ref{CRO}). Then, by combining the support property of the noise with Proposition \ref{CRO}, we complete the proof of irreducibility for the stochastic wave equation \eqref{WE}.

\subsection{Continuity of solving operator}
We recall that the space $\mathcal{G}^{2}$ is defined by
$\mathcal{G}^{2}:=\{\zeta\in C([0,T];H^{2}) \mid \zeta(0)=0\}.$
\begin{proposition}\label{CRO}
If $(u_{0},u_{1})\in \mathcal{H},g\in \mathcal{G}^{2},T>0,$ then the equation
\begin{eqnarray}\label{1}
 \begin{array}{l}
 \left\{
 \begin{array}{llll}
\ddot{u}+Au+\gamma \dot{u}+f(u)=\dot{g}
 \\u(x,0)=u_{0}(x)
 \\\dot{u}(x,0)=u_{1}(x)
 \end{array}
 \right.
 \end{array}
 \begin{array}{lll}
 {\rm{in}}~\mathbb{R}\times(0,T),\\
 {\rm{in}}~\mathbb{R},
 \\{\rm{in}}~\mathbb{R},
\end{array}
\end{eqnarray}
has a unique solution $(u,\dot{u})\in \mathcal{H}.$ Moreover, the map
\begin{equation*}
\begin{split}
F: \mathcal{G}^{2}&\mapsto C([0,T];\mathcal{H}),\\
g&\mapsto y^u=(u,\dot{u})
\end{split}
\end{equation*}
is well-defined and continuous.
\end{proposition}

\begin{proof}[Proof of Proposition \ref{CRO}]
Inspired by the framework developed in \cite{KS12,Gao22ECOCV,Gao26MA}, the proof is organized into the following steps.
\par
\textit{Step 1.} We prove the well-posedness of \eqref{1} and show that $F$ is locally bounded, that is , if $g\in \mathcal{G}^{2},\|g\|_{\mathcal{G}^{2}}\leq R,$
\begin{equation}\label{26}
\sup\limits_{0\leq t\leq T}\mathcal{E}_{u}(t)\leq C(T,R,u_{0},u_1).
\end{equation}
\par
Indeed, making the substitution $v(t):=u(t)-G(t)$, where $G(t)=\int_{0}^{t}g(s)ds$, we can see that
\begin{equation*}
\begin{split}
\ddot{v}+Av+\gamma \dot{v}=-f(v+G)-AG-\gamma \dot{G}:=h(v,G).
\end{split}
\end{equation*}
By virtue of the truncation technique as in the proof of Proposition \ref{WP}, we can prove the local well-posedness of \eqref{1}.
According to \eqref{E3}, we have
\begin{equation}\label{36}
\begin{split}
\frac{d}{dt}|y^{v}|_{\mathcal{H}}^{2}=-2\alpha\|v\|_{1}^{2}+2(\alpha-\gamma)\|\dot{v}+\alpha v\|^{2}-2\alpha(\alpha-\gamma)(v,\dot{v}+\alpha v)+\int_{\mathbb{R}}2(\dot{v}+\alpha v)h(v,G)dx.
\end{split}
\end{equation}
It follows from \eqref{59} that
$-f(v+G)=-|v|^{2m}v+D(v,G),$ where
$$|D(v,G)|\leq C_{m}(|v|^{2m}|G|+|G|^{2m+1}).$$
\eqref{36} implies that
\begin{equation*}
\begin{split}
&\frac{d}{dt}(|y^{v}|_{\mathcal{H}}^{2}+\frac{1}{m+1}\int_{\mathbb{R}}|v|^{2m+2}dx)\\
=&-2\alpha\|v\|_{1}^{2}+2(\alpha-\gamma)\|\dot{v}+\alpha v\|^{2}-2\alpha(\alpha-\gamma)(v,\dot{v}+\alpha v)
\\&-2\alpha\int_{\mathbb{R}}|v|^{2m+2}dx+\int_{\mathbb{R}}2(\dot{v}+\alpha v)D(v,G)dx
-\int_{\mathbb{R}}2(\dot{v}+\alpha v)[AG+\gamma \dot{G}]dx.
\end{split}
\end{equation*}
With the help of Young's inequality, we have
\begin{equation*}
\begin{split}
\int_{\mathbb{R}}2(\dot{v}+\alpha v)D(v,G)dx
&\leq \frac{\alpha}{2}\|\dot{v}+\alpha v\|^{2}+C_{\alpha,m}\int_{\mathbb{R}}(|v|^{4m}|G|^{2}+|G|^{4m+2})dx\\
&\leq \frac{\alpha}{2}\|\dot{v}+\alpha v\|^{2}+C_{\alpha,m}\int_{\mathbb{R}}(\delta|v|^{2m+2}+C_{\delta}|G|^{\frac{2(m+1)}{1-m}}+|G|^{4m+2})dx,\\
\int_{\mathbb{R}}2(\dot{v}+\alpha v)[AG+\gamma \dot{G}]dx
&\leq \frac{\alpha}{2}\|\dot{v}+\alpha v\|^{2}+C_{\alpha}\|AG+\gamma \dot{G}\|^{2},
\end{split}
\end{equation*}
this implies that
\begin{equation}\label{27}
\begin{split}
\frac{d}{dt}(|y^{v}|_{\mathcal{H}}^{2}+\frac{1}{m+1}\int_{\mathbb{R}}|v|^{2m+2}dx)
\leq&-\alpha(|y^{v}|_{\mathcal{H}}^{2}+\frac{1}{m+1}\int_{\mathbb{R}}|v|^{2m+2}dx)\\
&+C[\int_{\mathbb{R}}(|G|^{\frac{2(m+1)}{1-m}}+|G|^{4m+2})dx+\|AG(t)+\gamma \dot{G}(t)\|^{2}]
.
\end{split}
\end{equation}
Since $g\in \mathcal{G}^{2},$ applying Gronwall inequality to \eqref{27}, we can prove the global well-posedness of \eqref{1}, moreover, from the above analysis, we have \eqref{26}.
\par
\textit{Step 2.} We show that $F$ is continuous.
\par
Indeed, since $F(g_{i})=y^{u_{i}},$ we make the substitution $v_i(t):=u_i(t)-G_i(t)$, where $G_i(t):=\int_{0}^{t}g_i(s)ds, G:=G_{1}-G_{2}, g:=g_{1}-g_{2}$.
We define $v:=v_{1}-v_{2},$ it is easy to see that
\begin{equation*}
\ddot{v}+Av+\gamma \dot{v}=h(v_{1},G_{1})-h(v_{2},G_{2}).
\end{equation*}
According to \eqref{E3}, we have
\begin{equation*}
\begin{split}
\frac{d}{dt}|y^{v}|_{\mathcal{H}}^{2}=
&-2\alpha\|v\|_{1}^{2}+2(\alpha-\gamma)\|\dot{v}+\alpha v\|^{2}-2\alpha(\alpha-\gamma)(v,\dot{v}+\alpha v)\\
&+\int_{\mathbb{R}}2(\dot{v}+\alpha v)(h(v_{1},G_{1})-h(v_{2},G_{2}))dx.
\end{split}
\end{equation*}
With the help of \eqref{26}, it holds that
\begin{equation*}
\begin{split}
\|h(v_{1},G_{1})-h(v_{2},G_{2})\|\leq C(R,T)(\|v\|+\|g\|_{\mathcal{G}^{2}}).
\end{split}
\end{equation*}
Appying Gronwall's inequality, we have
\begin{equation*}
\begin{split}
\sup\limits_{0\leq t\leq T}|y^{v}|_{\mathcal{H}}\leq C(R,T)\|g\|_{\mathcal{G}^{2}}.
\end{split}
\end{equation*}
Since $u=v+G,$ we have
\begin{equation*}
\begin{split}
|y^{u_{1}}-y^{u_{2}}|_{\mathcal{H}}\leq C(u_{0},u_{1},T,R)\|g_{1}-g_{2}\|_{\mathcal{G}^{2}},
\end{split}
\end{equation*}
this implies that $F$ is continuous.
\par
This completes the proof.
\end{proof}

\subsection{Irreducibility for the stochastic wave equation}
Let us denote $\Phi$ as the \textit{resolving operator} for \eqref{1} as following
\begin{equation*}
\begin{split}
\Phi: \mathcal{H}\times \mathcal{G}^{2}&\rightarrow C([0,T];\mathcal{H})\\
(y_{0},g)&\mapsto y.
\end{split}
\end{equation*}
In what follows, we denote by $\Phi_{t}(y_{0},g)$ the restriction of $\Phi(y_{0},g)$ at time $t$.
That is, $\Phi_{t}$ takes $(y_{0},g)$ to $y(t)$, where $y(t)$ is the solution of \eqref{1}.
We consider system \eqref{WE} with $\eta$ in \eqref{37}. We define the \textit{transition probabilities} $\{P(t,y_{0},\cdot): t\in [0,T], y_{0}\in \mathcal{H}\}$ associated with system \eqref{WE} with $\eta$ in \eqref{37} by
$$P(t,y_{0},\Gamma):=\mathbb{P}\{\Phi_{t}(y_{0},ht+\int_{0}^{t}\eta(s)ds)\in \Gamma\}$$
for all Borel sets $\Gamma$ of $\mathcal{H}$ and all $t\in [0,T].$

\begin{proposition}\label{pro20}
Let the assumption $\textbf{(A)}$ hold. Then for any $R,d>0$, there exist constants $p_{0},T>0$ and $N\in \mathbb{N}$ such that if
\eqref{41} holds, we have
\begin{equation*}
\begin{split}
P(T,y_{0},B_{\mathcal{H}}(0,d))\geq p_{0}
\end{split}
\end{equation*}
for any $y_{0}\in B_{\mathcal{H}}(0,R)$.
\end{proposition}
\begin{proof}[Proof of Proposition \ref{pro20}]
Let $\bar{y}(t):=\Phi_{t}(y_{0},0)$. We also denote $\bar{y}=[\bar{u},\dot{\bar{u}}]$. If we take $f=-|\bar{u}|^{2m}\bar{u}$ in \eqref{E3}, then
\begin{equation}\label{64}
\begin{split}
\frac{d}{dt}\mathcal{E}_{\bar{u}}(t)
=-2\alpha\|\bar{u}\|_{1}^{2}+2(\alpha-\gamma)\|\dot{\bar{u}}+\alpha \bar{u}\|^{2}-2\alpha(\alpha-\gamma)(\bar{u},\dot{\bar{u}}+\alpha \bar{u})
-2\alpha\int_{\mathbb{R}}|\bar{u}|^{2m+2}dx.
\end{split}
\end{equation}
By taking $\alpha>0$ small enough, we have \eqref{25} holds. By substituting \eqref{25} into \eqref{64}, we arrive at
\begin{equation*}
\begin{split}
\frac{d}{dt}\mathcal{E}_{\bar{u}}(t)\leq-\frac{3}{2}\alpha|\bar{y}|_{\mathcal{H}}^{2}-2\alpha\int_{\mathbb{R}}|\bar{u}|^{2m+2}dx\leq-\frac{3}{2}\mathcal{E}_{\bar{u}}(t)
.
\end{split}
\end{equation*}
Gronwall's inequality then gives
\begin{equation}\label{65}
\begin{split}
\mathcal{E}_{\bar{u}}(t)\leq e^{-\frac{3}{2}\alpha t}\mathcal{E}_{\bar{u}}(0),~~{\rm{for}}~t\geq 0
.
\end{split}
\end{equation}
If $y_{0}\in B_{\mathcal{H}}(0,R)$, with the help of \eqref{65}, we can find a $T=T(R,d)$ such that \begin{equation}\label{39}
\begin{split}
|\bar{y}(T)|_{\mathcal{H}}< \frac{d}{2}.
\end{split}
\end{equation}
Applying Proposition~\ref{CRO} with $g(t):=th+W(t)$, we obtain the following: there exists $r>0$ such that, whenever
$\sup_{0\le t\le T}\|ht+W(t)\|_2 < r,$ we have
$|y(T)-\bar{y}(T)|_{\mathcal{H}} < \frac{d}{2}.$ This means that
\begin{equation*}
\begin{split}
\mathbb{P}(|y(T)-\bar{y}(T)|_{\mathcal{H}}< \frac{d}{2})\geq \mathbb{P}(\sup\limits_{0\leq t\leq T}\|th+W(t)\|_{2}< r).
\end{split}
\end{equation*}
If we choose $N$ sufficiently large, according to the assumption $\textbf{(A)}$ hold, we have
$$\|Q_{N}h\|_{2}\leq\sum\limits_{j\geq N+1}^{\infty}|(h,e_{j})|\|e_{j}\|_{2}<\frac{r}{3T},$$
this implies that
\begin{equation}\label{38}
\begin{split}
\mathbb{P}(|y(T)-\bar{y}(T)|_{\mathcal{H}}< \frac{d}{2})
\geq \mathbb{P}(\sup\limits_{0\leq t\leq T}\|P_{N}W(t)+tP_{N}h\|_{2}< \frac{r}{3})\cdot\mathbb{P}(\sup\limits_{0\leq t\leq T}\|Q_{N}W(t)\|_{2}< \frac{r}{3}).
\end{split}
\end{equation}
Assumption \eqref{41} implies that
$
p_1:=\mathbb{P}\left(\sup_{0\le t\le T}\|P_N W(t)+tP_N h\|_2 < \frac{r}{3}\right)>0.
$
On the other hand, since $\mathcal{B}_3<+\infty$, the distribution of $Q_N W$ is a centered Gaussian measure on $C([0,T];H^2)$. Hence,
$
p_2:=\mathbb{P}\left(\sup_{0\le t\le T}\|Q_N W(t)\|_2 < \frac{r}{3}\right)>0.
$
Combining these facts with \eqref{39} and \eqref{38}, we obtain
\begin{equation*}
\begin{split}
\mathbb{P}(|y(T)|_{\mathcal{H}}<d)\geq \mathbb{P}(|y(T)-\bar{y}(T)|_{\mathcal{H}}< \frac{d}{2})
\geq p_{1}p_{2}=:p_{0}.
\end{split}
\end{equation*}
This proves Proposition \ref{pro20}.
\end{proof}

\section{Proof of Theorem \ref{MT}}
By virtue of Theorem~\ref{Th1}, the proof of Theorem~\ref{MT} reduces to establishing recurrence and polynomial squeezing, and verifying condition (2) of Theorem~\ref{Th1}. The latter is already guaranteed by Proposition~\ref{WP}. It remains, therefore, to establish recurrence and polynomial squeezing.
To this end, we modify the framework and arguments of \cite{KS12,Shi08,M2014,N2,Gao26JDE} and adapt them to the stochastic wave equation on the whole line.

\subsection{Recurrence}
 We define the functional
\begin{equation*}
\begin{split}
V(y(t)):=\mathcal{E}_{u}(t)+1,
\end{split}
\end{equation*}
then, $V(y)$ is a continuous functional on $\mathcal{H}$ such that
\begin{equation*}
\begin{split}
V(y)\geq 1~{\rm{for~any}}~y\in \mathcal{H}~~{\rm{and}}~~\lim\limits_{|y|_{\mathcal{H}}\rightarrow+\infty}V(y)=+\infty.
\end{split}
\end{equation*}
Moreover, according to \eqref{45}, we can prove
\begin{equation*}
\begin{split}
&\mathbb{E}_{y}V(y(t_{*}))\leq aV(y),~~\forall |y|_{\mathcal{H}}\geq R,\\
&\mathbb{E}_{y}V(y(t))\leq C_{*},~~~~~~~~\forall |y|_{\mathcal{H}}\leq R,t\geq0,
\end{split}
\end{equation*}
for some constants $R>0,t_{*}>0,C_{*}>0,a<1.$ These mean that the family $(y(t), \mathbb{P}_{y})$ possesses
a Lyapunov function given by $V(y)$.
\par
Now, we state that the extension $(\mathbf{y}(t), \mathbb{P}_{\mathbf{y}})$ of $(y(t), \mathbb{P}_{y})$
is irreducible.
\begin{proposition}\label{pro11}
For any $R,d>0$, there exist constants $p,T>0$ and $N\in \mathbb{N}$ such that if
\eqref{41} holds, we have
\begin{equation}\label{34}
\begin{split}
\mathbb{P}_{\mathbf{y}}(\mathbf{y}(T)\in B_{\mathcal{H}}(0,d)\times B_{\mathcal{H}}(0,d))\geq p
\end{split}
\end{equation}
holds for any $\mathbf{y}=(y,y^\prime)\in B_{\mathcal{H}}(0,R)\times B_{\mathcal{H}}(0,R)$.
\end{proposition}

\begin{proof}[Proof of Proposition \ref{pro11}]
In order to prove Proposition~\ref{pro11}, we follow the framework in \cite{Shi08,M2014,N2,Gao26JDE} and adapt them to our setting.

Without loss of generality, we can assume $R\geq d.$
We define the events
\begin{equation*}
\begin{split}
&G_{d}(T):=\{|\mathcal{R}_{T}(y,y^\prime)|_{\mathcal{H}}\leq d\},\\
&G_{d}^{\prime}(T):=\{|\mathcal{R}_{T}^\prime(y,y^\prime)|_{\mathcal{H}}\leq d\},\\
&E_{\rho}:=\{\mathcal{F}_{\tilde{u}}^{\psi}(t)\leq M\mathcal{E}_{\tilde{u}}(0)+Kt+\rho\}
\cap \{\mathcal{F}_{\tilde{u}^\prime}^{\psi}(t)\leq M\mathcal{E}_{\tilde{u}^\prime}(0)+Kt+\rho\},\\
&\mathcal{N}:=\{y^{\tilde{v}}(t)\neq y^{\tilde{u}^{\prime}}(t)~{\rm{for~some~}}t\geq0\}.
\end{split}
\end{equation*}
We recall here the definitions of $\mathcal{R}_{T}$ and $\mathcal{R}_{T}^\prime$ introduced in Section~\ref{Sec2.2}.
\par
\textit{Step 1.}
We claim that
\begin{equation}\label{32}
\begin{split}
\mathbb{P}_{\mathbf{y}}(G_{ \frac{d}{2}}(kT))\wedge\mathbb{P}_{\mathbf{y}}(G^{\prime}_{ \frac{d}{2}}(kT))
\geq p_{0},
\end{split}
\end{equation}
for any $y\in B_{\mathcal{H}}(0,R)$ and $k\geq1.$
\par
In order to provide an accurate and detailed proof of inequality \eqref{32}, we introduce the notation $S_t(y)$ for the solution of \eqref{WE} issued from the initial data $y\in \mathcal{H}$.
Indeed, if \eqref{41} holds for a large enough integer $N\geq1$, it follows from Proposition \ref{pro20} that there exists a constant $p_{0}$ such that
$\mathbb{P}_{y}(|S_{T}(y)|_{\mathcal{H}}\leq \frac{d}{2})\geq p_{0}$ for any $y\in B_{\mathcal{H}}(0,R)$. This implies that
\begin{equation}\label{28}
\begin{split}
\mathbb{P}_{y}(|S_{kT}(y)|_{\mathcal{H}}\leq \frac{d}{2})\geq p_{0}
\end{split}
\end{equation}
for any $y\in B_{\mathcal{H}}(0,R)$ and $k\geq1.$ \eqref{28} means that \eqref{32} holds.
\par
Now, we prove \eqref{28} holds. Indeed, we introduce the stopping times
\begin{equation*}
\begin{split}
\bar{\tau}=\min\{nT:n\geq1,|S_{nT}(y)|_{\mathcal{H}}> \frac{d}{2}\},~\bar{\sigma}=\bar{\tau}\wedge (kT).
\end{split}
\end{equation*}
Since $\bar{\sigma}$ is a.s. finite, we can use the strong Markov property, and obtain
\begin{equation*}
\begin{split}
\mathbb{P}_{y}(|S_{kT}(y)|_{\mathcal{H}}> \frac{d}{2})
&\leq \mathbb{P}_{y}(\bar{\tau}\leq kT)=\mathbb{P}_{y}(\bar{\sigma}=\bar{\tau})
\leq \mathbb{P}_{y}(|S_{\bar{\sigma}-T}(y)|_{\mathcal{H}}\leq R, |S_{\bar{\sigma}}(y)|_{\mathcal{H}}> \frac{d}{2})\\
&=\mathbb{E}_{y}[\mathbb{E}_{y}(\mathbb{I}_{|S_{\bar{\sigma}-T}(y)|_{\mathcal{H}}\leq R}\cdot\mathbb{I}_{|S_{\bar{\sigma}}(y)|_{\mathcal{H}}\geq \frac{d}{2}}|\mathcal{F}_{\bar{\sigma}-T}) ]
=\mathbb{E}_{y}[\mathbb{I}_{|v|_{\mathcal{H}}\leq R}\cdot\mathbb{E}_{v}(\mathbb{I}_{|S_{T}(v)|_{\mathcal{H}}> \frac{d}{2}})]\\
&=\mathbb{E}_{y}[\mathbb{I}_{|v|_{\mathcal{H}}\leq R}\cdot\mathbb{P}_{v}(|S_{T}(v)|_{\mathcal{H}}> \frac{d}{2})]
\leq \sup_{\bar{v}\in B_{R}}\mathbb{P}_{\bar{v}}(|S_{T}(\bar{v})|_{\mathcal{H}}> \frac{d}{2})\\
&\leq 1-p_{0},
\end{split}
\end{equation*}
where $v=S_{\bar{\sigma}-T}(y)$.

\par
\textit{Step 2.} We claim that for any $\rho>0,$ there exists $k\geq 1$ such that
\begin{equation}\label{33}
\begin{split}
G_{\frac{d}{2}}(kT)E_{\rho}\mathcal{N}^{c}\subset G_{d}(kT)G_{d}^{\prime}(kT)
\end{split}
\end{equation}
for any $\mathbf{y}\in B_{\mathcal{H}}(0,R)\times B_{\mathcal{H}}(0,R).$
\par
Indeed, we only need to show that on the event $G_{\frac{d}{2}}(kT)E_{\rho}\mathcal{N}^{c},$
\begin{equation}\label{31}
\begin{split}
|\mathcal{R}_{kT}(y,y^\prime)-\mathcal{R}_{kT}^\prime(y,y^\prime)|_{\mathcal{H}}\leq\frac{d}{2}
\end{split}
\end{equation}
holds for large enough $k\geq 1$ and any $\mathbf{y}\in B_{\mathcal{H}}(0,R)\times B_{\mathcal{H}}(0,R).$
Since $\omega\in \mathcal{N}^{c},$ we have $\tilde{v}=\tilde{u}^{\prime},$ we have $\mathcal{R}_{t}^\prime(y,y^\prime)=\tilde{u}^{\prime}=\tilde{v}$ and $\mathcal{R}_{t}(y,y^\prime)=\tilde{u}$. Noting the fact $\tilde{v}$ and $\tilde{u}$ solve the equations \eqref{29} and \eqref{30}, their difference $\tilde{w}=\tilde{u}-\tilde{v}$ solves the equation \eqref{24} with $u,v$ replaced by $\tilde{u},\tilde{v}$. By virtue of Theorem~\ref{FP} and an argument parallel to that leading to \eqref{11}, it follows that
\begin{equation*}
\begin{split}
|\mathcal{R}_{kT}(y,y^\prime)-\mathcal{R}_{kT}^\prime(y,y^\prime)|_{\mathcal{H}}^{2}\leq C(\rho,R)e^{-\alpha kT}|y-y^{\prime}|_{\mathcal{H}}^{2}
.
\end{split}
\end{equation*}
Because $\rho$ is fixed, one can select $k \ge 1$ large enough so that the right-hand side of the inequality is bounded by $\frac{d^2}{4}$, which ensures that \eqref{31} is satisfied.
\par
\textit{Step 3. Proof of \eqref{34}.}
\par
Indeed, it follows from Proposition \ref{pro2} that
$
\mathbb{P}_{\mathbf{y}}(E_{\rho})\geq1-o(\rho).
$
Taking into account the conditional independence of $\tilde{v}$ and $\tilde{u}'$ given $\mathcal{N}$, and invoking \eqref{33}, we deduce that
\begin{equation}\label{35}
\begin{split}
\mathbb{P}_{\mathbf{y}}(G_{d}G_{d}^{\prime})
&=\mathbb{P}_{\mathbf{y}}(G_{d}G_{d}^{\prime}\mathcal{N}^{c})+\mathbb{P}_{\mathbf{y}}(G_{d}G_{d}^{\prime}\mathcal{N})\\
&\geq\mathbb{P}_{\mathbf{y}}(G_{\frac{d}{2}}E_{\rho}\mathcal{N}^{c})+\mathbb{P}_{\mathbf{y}}(G_{d}G_{d}^{\prime}\mathcal{N})\\
&\geq\mathbb{P}_{\mathbf{y}}(G_{\frac{d}{2}}\mathcal{N}^{c})+\mathbb{P}_{\mathbf{y}}(G_{d}\mathcal{N})\mathbb{P}_{\mathbf{y}}(G_{d}^{\prime}\mathcal{N})-o(\rho)
.
\end{split}
\end{equation}
By taking $\rho>0$ large enough such that $o(\rho)\leq \frac{p^{2}_{0}}{8}.$ By the symmetry, we can assume that $\mathbb{P}_{\mathbf{y}}(G_{\frac{d}{2}}^{\prime}\mathcal{N}^{c})\leq\mathbb{P}_{\mathbf{y}}(G_{\frac{d}{2}}\mathcal{N}^{c}).$ If $\mathbb{P}_{\mathbf{y}}(G_{\frac{d}{2}}\mathcal{N}^{c})\geq \frac{p^{2}_{0}}{4}$, \eqref{34} is proved.
If $\mathbb{P}_{\mathbf{y}}(G_{\frac{d}{2}}\mathcal{N}^{c})\leq \frac{p^{2}_{0}}{4}$, it follows from \eqref{32} that
\begin{equation*}
\begin{split}
p_{0}^{2}\leq\mathbb{P}_{\mathbf{y}}(G_{\frac{d}{2}})\mathbb{P}_{\mathbf{y}}(G_{\frac{d}{2}}^{\prime})
\leq\mathbb{P}_{\mathbf{y}}(G_{\frac{d}{2}}\mathcal{N})\mathbb{P}_{\mathbf{y}}(G_{\frac{d}{2}}^{\prime}\mathcal{N})+\frac{3}{4}p_{0}^{2},
\end{split}
\end{equation*}
this leads to
\begin{equation*}
\begin{split}
\mathbb{P}_{\mathbf{y}}(G_{d}\mathcal{N})\mathbb{P}_{\mathbf{y}}(G_{d}^{\prime}\mathcal{N})\geq \mathbb{P}_{\mathbf{y}}(G_{\frac{d}{2}}\mathcal{N})\mathbb{P}_{\mathbf{y}}(G_{\frac{d}{2}}^{\prime}\mathcal{N})\geq \frac{1}{4}p_{0}^{2},
\end{split}
\end{equation*}
by putting this estimate into \eqref{35}, we can prove \eqref{34}.
\par
This completes the proof.
\end{proof}

It follows from \cite[Proposition 3.3]{Shi08} that the above facts imply that recurrence in Theorem \ref{Th1} holds for the extension $(\mathbf{y}(t),\mathbb{P}_{\mathbf{y}}).$ Summarizing the above, we arrive at the following result.
\begin{proposition}\label{Rec}
For any $d>0$, there exist constants $c,C>0$ and $N\in \mathbb{N}$ such that if
\eqref{41} holds, we have
\begin{equation*}
\begin{split}
\mathbb{E}_{\mathbf{y}}e^{c\tau_{\textbf{B}}}\leq C(V(y)+V(y^\prime))~~{\rm{for~any~}}\mathbf{y}=(y,y^\prime)\in\mathcal{H}\times\mathcal{H}.
\end{split}
\end{equation*}
where $B:=B_{\mathcal{H}}(0,d)$. In particular, recurrence in Theorem~\ref{Th1} holds for the extension $(\mathbf{y}(t),\mathbb{P}_{\mathbf{y}})$.
\end{proposition}

\subsection{Polynomial squeezing}
Following the nations in \cite{KS12,Shi08,M2014,N2,Gao26JDE}, let us define the stopping time $\sigma$ as follows
\begin{equation*}
\begin{split}
\sigma:=\tilde{\tau}\wedge \theta,
\end{split}
\end{equation*}
where $\tilde{\tau}:=\tau^{\tilde{u}}\wedge \tau^{\tilde{u}^{\prime}}$ and $\theta:=\inf\{t\geq0~|~y^{\tilde{u}^{\prime}}(t)\neq y^{\tilde{v}}(t)\}$
and introduce the following events
\begin{equation*}
\begin{split}
\mathcal{Q}_{k}^{\prime}:=\{kT\leq \sigma\leq(k+1)T:\theta\geq \tilde{\tau}\},\\
\mathcal{Q}_{k}^{\prime\prime}:=\{kT\leq \sigma\leq(k+1)T:\theta<\tilde{\tau}\}.
\end{split}
\end{equation*}
For $\mathcal{Q}_{k}^{\prime},\mathcal{Q}_{k}^{\prime\prime},$ we have the following property.
\begin{proposition}\label{pro12}
For any $q>1,$ there exist $N\in \mathbb{N}$, large enough constants $\rho,L,T>0$ and small enough $d>0$ such that
\begin{equation*}
\begin{split}
\mathbb{P}_{\mathbf{y}}(\mathcal{Q}_{k}^{\prime})\vee\mathbb{P}_{\mathbf{y}}(\mathcal{Q}_{k}^{\prime\prime}) \leq \frac{q-1}{4q}\frac{1}{(k+2)^{q}}
\end{split}
\end{equation*}
for any $\mathbf{y}\in \bar{B}_{\mathcal{H}}(0,d)\times\bar{B}_{\mathcal{H}}(0,d),$ provided that \eqref{41} holds.
\end{proposition}

\begin{proof}[Proof of Proposition \ref{pro12}]
In order to prove Proposition~\ref{pro12}, we modify the framework and arguments of \cite{KS12,Shi08,M2014,N2,Gao26JDE} and adapt them to our setting.

The proof is divided into several steps.
\par
\textit{Step 1. Estimate of $\mathbb{P}_{\mathbf{y}}(\mathcal{Q}_{k}^{\prime})$.}
\par
By choosing $\rho$ and $L$ sufficiently large and applying Proposition~\ref{pro1} with $q$ replaced by $q+1$, we obtain
\begin{equation*}
\begin{split}
\mathbb{P}_{\mathbf{y}}(\mathcal{Q}_{k}^{\prime})\leq \mathbb{P}_{\mathbf{y}}(kT\leq \tilde{\tau}<+\infty)
\leq \mathbb{P}_{y}(kT\leq \tau^{\tilde{u}}<+\infty)\leq 2Q_{q+1}(\rho+kLT+1)\leq \frac{q-1}{4q}\frac{1}{(k+2)^{q}},
\end{split}
\end{equation*}

\par
\textit{Step 2. Estimate of $\mathbb{P}_{\mathbf{y}}(\mathcal{Q}_{0}^{\prime\prime})$.}
\par
For any $\mathbf{y}=(y,y^{\prime})\in \bar{B}_{\mathcal{H}}(0,d)\times\bar{B}_{\mathcal{H}}(0,d),$ by the definition, we have
\begin{equation*}
\begin{split}
\mathbb{P}_{\mathbf{y}}(\mathcal{Q}_{0}^{\prime\prime})
\leq \mathbb{P}_{\mathbf{y}}\{0\leq \theta \leq T\}=\mathbb{P}_{\mathbf{y}}\{y^{\tilde{u}^{\prime}}(t)\neq y^{\tilde{v}}(t)~{\rm{for~some}}~t\in[0,T]\}
=\|\lambda_T(y,y^{\prime})-\lambda_T^{\prime}(y,y^{\prime})\|_{var},
\end{split}
\end{equation*}
where $\lambda_T(y,y')$ and $\lambda'_T(y,y')$ denote the distributions of $\{y^{v}(t)\}_{t\in[0,T]}$ and $\{y^{u'}(t)\}_{t\in[0,T]}$ respectively. $(V_{T}(y,y^{\prime}),V_{T}^{\prime}(y,y^{\prime}))$ is a maximal coupling for the pair $(\lambda_{T}(y,y^{\prime}),\lambda_{T}^{\prime}(y,y^{\prime}))$, and we use the fact that $\{y^{\tilde{v}}(t)\}_{t\in[0,T]}$ and $\{y^{\tilde{u}'}(t)\}_{t\in[0,T]}$ are flows of the maximal coupling $(V_T(y,y'),V_T'(y,y'))$. We define $$\tau:=\tau^{u}\wedge\tau^{u'}\wedge\tau^v.$$
Then, we can see that
\begin{equation*}
\begin{split}
\|\lambda_T(y,y^{\prime})-\lambda_T^{\prime}(y,y^{\prime})\|_{var}
&=\sup_{\Gamma\in \mathcal{B}(C([0,T];\mathcal{H}))}|\mathbb{P}\{y^{v}(\cdot)\in \Gamma\}-\mathbb{P}\{y^{u'}(\cdot)\in \Gamma\}|\\
&\leq \mathbb{P}(\tau<\infty)
+\sup_{\Gamma\in \mathcal{B}(C([0,T];\mathcal{H}))}|\mathbb{P}\{y^{v}(\cdot)\in \Gamma, \tau=\infty\}-\mathbb{P}\{y^{u'}(\cdot)\in \Gamma, \tau=\infty\}|
.
\end{split}
\end{equation*}
It is not difficult to obtain that
\begin{equation*}
\begin{split}
\mathbb{P}(\tau<\infty)
\leq \mathbb{P}\{\tau^{u}<\infty\}+\mathbb{P}\{\tau^{u'}<\infty\}+\mathbb{P}\{\tau^{v}<\infty\}.
\end{split}
\end{equation*}
According to the definitions of $\hat{u}'$ and $\hat{v}$, it holds that
\begin{equation*}
\begin{split}
&\sup_{\Gamma\in \mathcal{B}(C([0,T];\mathcal{H}))}|\mathbb{P}\{y^{v}(\cdot)\in \Gamma, \tau=\infty\}-\mathbb{P}\{y^{u'}(\cdot)\in \Gamma, \tau=\infty\}|
\\=&\sup_{\Gamma\in \mathcal{B}(C([0,T];\mathcal{H}))}|\mathbb{P}\{y^{\hat{v}}(\cdot)\in \Gamma, \tau=\infty\}-\mathbb{P}\{y^{\hat{u}'}(\cdot)\in \Gamma, \tau=\infty\}|
\\\leq&\sup_{\Gamma\in \mathcal{B}(C([0,T];\mathcal{H}))}|\mathbb{P}\{y^{\hat{v}}(\cdot)\in \Gamma\}-\mathbb{P}\{y^{\hat{u}'}(\cdot)\in \Gamma\}|
\\\leq&\|\mathbb{P}-\Phi^{u,u^{\prime}}_{\ast}\mathbb{P}\|_{var},
\end{split}
\end{equation*}
where $\Phi^{u,u'}$ is the transformation defined in \eqref{58}. According to the above analysis, we have
\begin{equation}\label{23}
\begin{split}
\mathbb{P}_{\mathbf{y}}(\mathcal{Q}_{0}^{\prime\prime})
&\leq\mathbb{P}_{\mathbf{y}}(0\leq \theta\leq T)\\
&\leq 2\mathbb{P}_{\mathbf{y}}(\tau^{u}<\infty)
+2\mathbb{P}_{\mathbf{y}}(\tau^{u^{\prime}}<\infty)
+\mathbb{P}_{\mathbf{y}}(\tau^{\hat{u}^{\prime}}<\infty)
+2\|\mathbb{P}-\Phi^{u,u^{\prime}}_{\ast}\mathbb{P}\|_{var}
.
\end{split}
\end{equation}
From \eqref{23}, together with Proposition~\ref{pro0} and Proposition~\ref{pro01} applied with $q+1$ instead of $q$, and using \eqref{13}, we deduce
\begin{equation*}
\begin{split}
\mathbb{P}_{\mathbf{y}}(\mathcal{Q}_{0}^{\prime\prime})
\leq CQ_{q+1}(\rho+1)+C\left(\exp(Ce^{C(\rho+(d^{2}+d^{2m+2})^2)}d^{2})-1 \right)^{\frac{1}{2}}.
\end{split}
\end{equation*}
By taking $\rho$ large enough and $d$ small enough, we have
\begin{equation*}
\begin{split}
\mathbb{P}_{\mathbf{y}}(\mathcal{Q}_{0}^{\prime\prime})\leq \frac{q-1}{4q}\frac{1}{2^q}.
\end{split}
\end{equation*}
\par
\textit{Step 3. Estimate of $\mathbb{P}_{\mathbf{y}}(\mathcal{Q}_{k}^{\prime\prime})$.}
\par
With the help of the Markov property and the estimate for $\mathbb{P}_{\mathbf{y}}(\mathcal{Q}_{0}^{\prime\prime})$,
we establish the estimate of $\mathbb{P}_{\mathbf{y}}(\mathcal{Q}_{k}^{\prime\prime})(k\geq1)$. By the Markov property we have
\begin{equation*}
\begin{split}
\mathbb{P}_{\mathbf{y}}(\mathcal{Q}_{k}^{\prime\prime})
=\mathbb{P}_{\mathbf{y}}(\mathcal{Q}_{k}^{\prime\prime},\sigma\geq kT)
=\mathbb{E}_{\mathbf{y}}(\mathbb{I}_{\sigma\geq kT}\cdot\mathbb{E}_{\mathbf{y}}(\mathbb{I}_{\mathcal{Q}_{k}^{\prime\prime}}|\mathcal{F}_{kT}))
\leq\mathbb{E}_{\mathbf{y}}(\mathbb{I}_{\sigma\geq kT}\cdot \mathbb{P}_{\mathbf{y}(kT)}(0\leq \theta\leq T))
,
\end{split}
\end{equation*}
combining this with \eqref{23}, we have
\begin{equation*}
\begin{split}
\mathbb{P}_{\mathbf{y}}(\mathcal{Q}_{k}^{\prime\prime})
\leq &2\mathbb{E}_{\mathbf{y}}(\mathbb{I}_{\sigma\geq kT}\cdot \mathbb{P}_{\mathbf{y}(kT)}(\tau^{u}<\infty))
+2\mathbb{E}_{\mathbf{y}}(\mathbb{I}_{\sigma\geq kT}\cdot \mathbb{P}_{\mathbf{y}(kT)}(\tau^{u^{\prime}}<\infty))\\
&+\mathbb{E}_{\mathbf{y}}(\mathbb{I}_{\sigma\geq kT}\cdot \mathbb{P}_{\mathbf{y}(kT)}(\tau^{\hat{u}^{\prime}}<\infty))
+2\mathbb{E}_{\mathbf{y}}(\mathbb{I}_{\sigma\geq kT}\cdot \|\mathbb{P}-\Phi^{\tilde{y}(kT),\tilde{y}^{\prime}(kT)}_{\ast}\mathbb{P}\|_{var})\\
:=&2I_{1}+2I_{2}+I_{3}+2I_{4},
\end{split}
\end{equation*}
where
\begin{equation*}
\begin{split}
\Phi^{\tilde{y}(kT),\tilde{y}'(kT)}(\omega)_t=\omega_t-\int_0^t\mathbb{I}_{\{s\le \tau^{k}\}}\mathcal{P}_N (0,f(u_k)-f(v_k))ds,
\end{split}
\end{equation*}
and $y_{k}=[u_{k},\dot{u}_{k}],y_{k}^{\prime}=[u^{\prime}_{k},\dot{u}^{\prime}_{k}]$ are the solutions of \eqref{WE} issued from $\tilde{y}(kT),\tilde{y}^{\prime}(kT),$ respectively.

\par
Let us estimate $I_{i}(i=1,2,3,4).$
\par
For the term $I_{4},$ let $z_{k}=[v_{k},\dot{v}_{k}]$ be the solution of \eqref{F1} issued from $\tilde{y}^{\prime}(kT),$ we define
$$\tau^{k}:=\tau^{u_{k}}\wedge\tau^{u_{k}^{\prime}}\wedge\tau^{v_{k}}.$$
In order to bound $I_{4}$, we need to estimate the integral
$$
\mathcal{A}_k(t):=-\mathbb{I}_{\{s\le \tau^{k}\}}\mathcal{P}_N (0,f(u_k)-f(v_k))ds.
$$
To this end, let us derive a pathwise estimate for $w_k(t):=u_k(t)-v_k(t)$. Let $T_{0}$ be in Theorem \ref{FP}.
\par
\textit{Case 1. $\tau^{k}\leq T_{0}$.}
\par
By applying (1) in Theorem \ref{FP}, we have
\begin{equation}\label{52}
\begin{split}
|y^{w_{k}}(t)|_{\mathcal{H}}^{2}
\leq |y^{w_{k}}(0)|_{\mathcal{H}}^{2}\exp\left (-\alpha t+C\int_{0}^{t}(\|u_{k}\|_{1}^{2}+\|v_{k}\|_{1}^{2})ds\right)
.
\end{split}
\end{equation}
For $t\leq \tau^{k},$ by the definition of $\tau^{k}$, we have
\begin{equation}\label{53}
\begin{split}
\mathcal{F}_{u_{k}}^{\psi}(t)\leq M\mathcal{E}_{u_{k}}(0)+(K+L)t+\rho,~~\mathcal{F}_{v_{k}}^{\psi}(t)\leq M\mathcal{E}_{v_{k}}(0)+(K+L)t+\rho.
\end{split}
\end{equation}
It follows from \eqref{52} and \eqref{53} that
\begin{equation*}
\begin{split}
|y^{w_{k}}(t)|_{\mathcal{H}}^{2}
\leq &|y^{w_{k}}(0)|_{\mathcal{H}}^{2}\exp\left (-\alpha t+C[M(\mathcal{E}_{u_{k}}(0)+\mathcal{E}_{v_{k}}(0))+(K+L)T_{0}+\rho]\right)\\
\leq &C|y^{w_{k}}(0)|_{\mathcal{H}}^{2}\exp\left (-\alpha t+C(\mathcal{E}_{u_{k}}(0)+\mathcal{E}_{v_{k}}(0))\right)
.
\end{split}
\end{equation*}
With the help of Young inequality, we have
\begin{equation}\label{54}
\begin{split}
|y^{w_{k}}(t)|_{\mathcal{H}}^{2}
\leq &C|y^{w_{k}}(0)|_{\mathcal{H}}^{2}\exp\left (-\alpha t+\frac{\alpha}{100(K+L)}(\mathcal{E}_{u_{k}}^{2}(0)+\mathcal{E}_{v_{k}}^{2}(0))\right)\\
=&C|y^{w_{k}}(0)|_{\mathcal{H}}^{2}\exp\left (-\alpha t+\frac{\alpha}{100(K+L)}(\mathcal{E}_{\tilde{u}}^{2}(kT)+\mathcal{E}_{\tilde{v}}^{2}(kT))\right)
.
\end{split}
\end{equation}
\par
\textit{Case 2. $\tau^{k}> T_{0}$.}
\par
It is clear \eqref{54} holds for $[0,T_{0}]$. By applying (2) in Theorem \ref{FP} with $\varepsilon=\frac{\alpha}{100C_{\ast} (K+L)}(1\wedge \frac{1}{2M})$, it holds that
\begin{equation*}
\begin{split}
|y^{w_{k}}(t)|_{\mathcal{H}}^{2}
&\leq C|y^{w_{k}}(T_{0})|_{\mathcal{H}}^{2}\exp\left (-\frac{\alpha}{2}t+\frac{\alpha}{100(K+L)}(\mathcal{E}_{u_{k}}^{2}(0)+\mathcal{E}_{v_{k}}^{2}(0))\right)\\
&\leq C|y^{w_{k}}(0)|_{\mathcal{H}}^{2}\exp\left (-\frac{\alpha}{2} t+\frac{\alpha}{50(K+L)}(\mathcal{E}_{\tilde{u}}^{2}(kT)+\mathcal{E}_{\tilde{v}}^{2}(kT))\right)
.
\end{split}
\end{equation*}
On the set $\{\sigma\ge kT\}$, it holds that
\begin{equation*}
\begin{split}
&\mathcal{E}_{\tilde{u}}^{2}(kT)\leq\mathcal{F}_{u_{k}}^{\psi,2}(kT)< (K+L)kT+\rho+M\mathcal{E}_{\tilde{u}}^{2}(0),\\
&\mathcal{E}_{\tilde{v}}^{2}(kT)\leq\mathcal{F}_{v_{k}}^{\psi,2}(kT)< (K+L)kT+\rho+M\mathcal{E}_{\tilde{v}}^{2}(0),\\
&|y^{w_{k}}(0)|_{\mathcal{H}}^{2}= |y^{\tilde u}(kT)-y^{\tilde u'}(kT)|_{\mathcal{H}}^{2}\le C d^2 e^{-\frac{\alpha kT}{2}}e^{C(\rho+(d^{2}+d^{2m+2})^{2})}.
\end{split}
\end{equation*}
Combining the above results together, we conclude that on the set $\{\sigma\ge kT\}$, for $t\leq \tau^{k},$
\begin{equation*}
\begin{split}
|y^{w_{k}}(t)|_{\mathcal{H}}^{2}\leq Ce^{-\frac{\alpha t}{2}}e^{-\frac{\alpha kT}{4}}e^{C(\rho+d^{4m+4})}d^2.
\end{split}
\end{equation*}
Applying the similar argument as in Proposition \ref{pro10}, we have
\begin{equation*}
\begin{split}
\int_0^{\infty}\|\mathcal{A}_k(t)\|^2 dt\leq Ce^{C(\rho+d^{4m+4})}d^2e^{-\frac{\alpha kT}{8}}
\end{split}
\end{equation*}
on the set $\{\sigma\ge kT\}$.
Hence,
\begin{equation*}
\begin{split}
I_4\le \left(\exp\left(Ce^{C(\rho+d^{4m+4})}d^2e^{-\frac{\alpha kT}{8}}\right)-1\right)^{\frac{1}{2}}.
\end{split}
\end{equation*}
By taking $d<1$ and $T$ large enough, we have
\begin{equation*}
\begin{split}
I_4 \leq \frac{q-1}{32q}\frac{1}{(k+2)^{q}},\quad  k\ge 1.
\end{split}
\end{equation*}
\par
For the term $I_1$, we apply Proposition~\ref{pro1} with $q+1$ instead of $q$ and invoke the Markov property to obtain
\begin{equation*}
\begin{split}
\frac{q-1}{32q}\frac{1}{(k+2)^{q}}
&\geq\mathbb{P}_{\mathbf{y}}(kT\leq \tau^{u}<\infty)
=\mathbb{E}_{\mathbf{y}}[\mathbb{E}_{\mathbf{y}}(\mathbb{I}_{kT\leq \tau^{u}<\infty}|\mathcal{F}_{kT})]
=\mathbb{E}_{\mathbf{y}}[\mathbb{P}_{\mathbf{y}(kT)}(\tau^{u}<\infty)]\\
&\geq\mathbb{E}_{\mathbf{y}}[\mathbb{I}_{\sigma\geq kT}\cdot\mathbb{P}_{\mathbf{y}(kT)}(\tau^{u}<\infty)]
=I_{1}.
\end{split}
\end{equation*}
\par
For the term $I_{2},$ we also have
$
I_{2}\leq \frac{q-1}{32q}\frac{1}{(k+2)^{q}}.
$ According to Proposition \ref{proT2}, it holds that
$
I_{3}\leq \frac{q-1}{16q}\frac{1}{(k+2)^{q}}.
$
Gathering together the above estimates for $I_{i}$, we conclude
\begin{equation*}
\begin{split}
\mathbb{P}_{\mathbf{y}}(\mathcal{Q}_{k}^{\prime\prime})
\leq \frac{q-1}{4q}\frac{1}{(k+2)^{q}}.
\end{split}
\end{equation*}
This completes the proof.
\end{proof}

\par
We now prove that the polynomial squeezing property in Theorem~\ref{Th1} holds for the extension $(\mathbf{y}(t),\mathbb{P}_{\mathbf{y}})$.
\par
For any $p>1$, we take $q^{\prime}>p+1.$
We fix $N$ large enough so that Proposition~\ref{pro12}, Proposition~\ref{Rec} and (2) of Theorem~\ref{FP} hold with $\varepsilon=\frac{\alpha}{100C_{\ast}(K+L)}(1\wedge \frac{1}{M})$. Then, by virtue of the definition of $\sigma$ and condition (2) of Theorem~\ref{FP}, we can prove that the first item in the polynomial squeezing property of Theorem~\ref{Th1} holds.
Proposition \ref{pro12} (with $q=q^{\prime}$) implies that
\begin{equation*}
\begin{split}
\mathbb{P}_{\mathbf{y}}(\sigma=\infty)
\geq 1-\sum_{k=0}^{\infty}\mathbb{P}_{\mathbf{y}}(\sigma\in [kT,(k+1)T])
\geq 1-2\sum_{k=0}^{\infty}\frac{q^{\prime}-1}{4q^{\prime}}\frac{1}{(k+2)^{q^{\prime}}}:=c_{2}\geq \frac{1}{2}.
\end{split}
\end{equation*}
It follows from Proposition \ref{pro12} that
\begin{equation}\label{42}
\begin{split}
\mathbb{E}_{\mathbf{y}}(\mathbb{I}_{\sigma<\infty}\sigma^{p})
\leq \sum_{k=0}^{\infty}\mathbb{E}_{\mathbf{y}}(\mathbb{I}_{\sigma\in [kT,(k+1)T]}\sigma^{p})
\leq \sum_{k=0}^{\infty}\frac{q^{\prime}-1}{4q^{\prime}}\frac{1}{(k+2)^{q^{\prime}}}(k+1)^{p}T^{p}:=c_3.
\end{split}
\end{equation}
We define $g(|y|_{\mathcal{H}}):=V(y)$ and $G(\mathbf{y}):=g(|y^{\tilde{u}}|_{\mathcal{H}})+g(|y^{\tilde{u}^\prime}|_{\mathcal{H}})$.
From the definitions of $\sigma$, $\tau_1^{\tilde{u}}$ and $\tau_1^{\tilde{u}^{\prime}}$, if $\sigma<\infty$, we have
\begin{equation*}
\begin{split}
&\mathcal{E}_{\tilde{u}}(\sigma)\leq\mathcal{F}_{\tilde{u}}^{\psi}(\sigma)< (K+L)\sigma+\rho+M\mathcal{E}_{\tilde{u}}(0),\\
&\mathcal{E}_{\tilde{u}^{\prime}}(\sigma)\leq\mathcal{F}_{\tilde{u}^{\prime}}^{\psi}(\sigma)< (K+L)\sigma+\rho+M\mathcal{E}_{\tilde{u}^{\prime}}(0),
\end{split}
\end{equation*}
thus, for $1\leq q\leq p,$ we have
\begin{equation*}
\begin{split}
\mathbb{E}_{\mathbf{y}}(\mathbb{I}_{\sigma<\infty}G(\mathbf{y}(\sigma))^{q})
\leq  C\mathbb{E}_{\mathbf{y}}(\mathbb{I}_{\sigma<\infty}(1+\sigma)^{q})
\leq  C\mathbb{E}_{\mathbf{y}}(\mathbb{I}_{\sigma<\infty}\sigma^{p})
\leq  c_4,
\end{split}
\end{equation*}
where we have used \eqref{42}.

This completes the proof.
\section{Appendix}
\begin{proof}[Proof of Proposition \ref{WP}]
In order to prove Proposition~\ref{WP}, we borrow some ideas from \cite{Chow1,Chow2,M2014} and adapt them to our setting.
The proof is divided into the following steps.
\par
\par
\textit{Step 1.} For any $y_0\in\mathcal{H},$ the system \eqref{WE} has a unique continuous local solution $y\in\mathcal{H}$.
\par
For any $R>0,$ let $\rho\in C^{\infty}_{0}(\mathbb{R})$ be a cut-off function such that
\begin{eqnarray*}
 \begin{array}{l}
\eta_{N}(s)=
\left\{
 \begin{array}{llll}
1,\\
\in(0,1]\\
0,
 \end{array}
 \right.
 \end{array}
 \begin{array}{lll}
0\leq s\leq\frac{N}{2},\\
\frac{N}{2}\leq s\leq N,\\
s\geq N.
\end{array}
\end{eqnarray*}
We define $f_{N}(u):=\eta_{N}(\|u\|_{1})f(\eta_{N}(\|u\|_{1})u)$. With the help of the truncated term $f_{N}(u)$, we consider the truncated system defined below
\begin{eqnarray}\label{TWE}
 \begin{array}{l}
 \left\{
 \begin{array}{llll}
\ddot{u}+Au+\gamma \dot{u}+f_{N}(u)=h+\eta
 \\u(x,0)=u_{0}(x)
 \\\dot{u}(x,0)=u_{1}(x)
 \end{array}
 \right.
 \end{array}
 \begin{array}{lll}
 {\rm{in}}~\mathbb{R}\times(0,+\infty),\\
 {\rm{in}}~\mathbb{R},
 \\{\rm{in}}~\mathbb{R},
\end{array}
\end{eqnarray}
we can rewrite \eqref{TWE} as
\begin{eqnarray}\label{TWEY}
 \begin{array}{l}
 \left\{
 \begin{array}{llll}
dy=[\mathbb{A}y+\mathbb{F}_{N}(y)]dt+\mathbb{G}dW
 \\y(0)=y_{0}
 \end{array}
 \right.
 \end{array}
 \begin{array}{lll}
 {\rm{in}}~\mathbb{R}\times(0,+\infty),\\
 {\rm{in}}~\mathbb{R},
\end{array}
\end{eqnarray}
where $\mathbb{F}_N(y)=[0,-f_N(u)+h].$ According to the definition of cut-off function $\eta_{N}$, we have
\begin{equation*}
\begin{split}
&\|f_N(u_1)-f_N(u_2)\|\leq C(N)\|u_1-u_2\|_{1},\\
&\|f_N(u)\|\leq C(N)\|u\|_{1},
\end{split}
\end{equation*}
this implies that
\begin{equation*}
\begin{split}
&|\mathbb{F}_N(y_1)-\mathbb{F}_N(y_2)|_{\mathcal{H}}\leq C(N)|y_1-y_2|_{\mathcal{H}},\\
&|\mathbb{F}_N(y)|_{\mathcal{H}}\leq C(N,h)(|y|_{\mathcal{H}}+1).
\end{split}
\end{equation*}
Since $\mathbb{F}_N$ are uniformly Lipschitz continuous and of linear growth, then
with the help of the standard existence and uniqueness theorem (see \cite{Chow2015,DaZ1}) for a stochastic evolution equation, we can obtain that
there exists a unique continuous solution $y_{N}(t):=[u_{N}(t),\dot{u}_{N}(t)]$ to the
truncated system \eqref{TWE}. To construct a local solution, we introduce the stopping time
\[
\tau_N := \inf \left\{ t > 0 : \|u_N(t)\|_1 > \frac{N}{2} \right\}.
\]
On the interval $t < \tau_N$, the function $u_N(t)$ coincides with the solution $u(t)$ of \eqref{WE}. Since $\tau_N$ is increasing with respect to $N$, the limit
\[
\tau_\infty := \lim_{N \to \infty} \tau_N
\]
is well defined. We then define $u(t)$ for $t < \tau_\infty \wedge T$ by setting $u(t) = u_N(t)$ whenever $t < \tau_N \le T$. This defines a unique local solution of \eqref{WE}.
\par
\textit{Step 2.} For any fixed $T > 0$ and initial datum $y_0 \in \mathcal H$, the Cauchy problem \eqref{WE} has a unique continuous solution $y(t) \in \mathcal H$ defined on $[0,T]$. We will demonstrate that the local solution obtained above is, in fact, global. To do so, it suffices to prove that the stopping times $\tau_N$ tend to infinity almost surely as $N \to \infty$. Indeed, once this is established, we have $u(t \wedge \tau_N) \to u(t)$ a.s. for every $t \le T$, which implies that the solution exists for all $t \ge 0$.
\par
Indeed, recall that on the time interval $t < \tau_N \le T$, $u(t)$ coincides with $u_N(t)$, and the latter satisfies system \eqref{WE}. It follows from \eqref{19} that
\begin{equation}\label{43}
\begin{split}
\mathcal{E}_{u}(t\wedge \tau_{N})+\frac{3}{2}\int_{0}^{t\wedge \tau_{N}}\alpha\mathcal{E}_{u}(s)ds
\leq M(t\wedge \tau_{N})+(\|h\|^{2}+\mathcal{B}_{1})t
.
\end{split}
\end{equation}
By taking the expectation of the above inequality, we have
\begin{equation*}
\begin{split}
\mathbb{E}\mathcal{E}_{u}(t\wedge \tau_{N})+\frac{3\alpha}{2}\mathbb{E}\int_{0}^{t\wedge \tau_{N}}\mathcal{E}_{u}(s)ds
\leq (\|h\|^{2}+\mathcal{B}_{1})t
,
\end{split}
\end{equation*}
which implies that for any $T>0$, we have
\begin{equation*}
\begin{split}
\mathbb{E}\mathcal{E}_{u}(T\wedge \tau_{N})\geq\mathbb{E}\mathbb{I}_{(\tau_{N}\leq T)}\mathcal{E}_{u}(\tau_{N})
\geq\mathbb{E}\mathbb{I}_{(\tau_{N}\leq T)}\|u(\tau_{N})\|_{1}^{2}\geq (\frac{N}{2})^{2}\mathbb{P}(\tau_{N}\leq T)
,
\end{split}
\end{equation*}
this leads to
\begin{equation*}
\begin{split}
\mathbb{P}(\tau_{N}\leq T)\leq \frac{4\mathbb{E}\mathcal{E}_{u}(T\wedge \tau_{N})}{N^{2}}
\leq \frac{C(T)}{N^{2}}
.
\end{split}
\end{equation*}
By virtue of the Borel-Cantelli lemma, this implies that
$
\mathbb{P}(\tau_{\infty}\leq T)=0
,
$
namely, we have
$
\lim\limits_{N\rightarrow\infty}\tau_{N}=\infty~~{\rm{a.s.}}
$
We now construct the global solution by setting
\[
u(t) := \lim_{N \to \infty} u(\tau_N \wedge t),
\]
which is well defined thanks to the previous convergence result. This $u(t)$ is the global solution we seek.
\par
Finally, we prove \eqref{45}.
\par
Indeed, it follows from \eqref{E1} that
\begin{equation*}
\begin{split}
d\mathcal{E}_{u}(t)=&[-2\alpha\|u\|_{1}^{2}+2(\alpha-\gamma)\|\dot{u}+\alpha u\|^{2}-2\alpha(\alpha-\gamma)(u,\dot{u}+\alpha u)-2\alpha\int_{\mathbb{R}}|u|^{2m+2}dx+\|h\|^{2}+\mathcal{B}_{1}]dt\\
&+2\sum\limits_{j=1}^{\infty}b_{j}(\dot{u}+\alpha u,e_{j})d\beta_{j}(t),
\end{split}
\end{equation*}
this gives that
\begin{equation*}
\begin{split}
\mathcal{E}_{u}(t\wedge \tau_{N})=&\mathcal{E}_{u}(0)+\int_{0}^{t\wedge \tau_{N}}[-2\alpha\|u\|_{1}^{2}+2(\alpha-\gamma)\|\dot{u}+\alpha u\|^{2}-2\alpha(\alpha-\gamma)(u,\dot{u}+\alpha u)\\
&-2\alpha\int_{\mathbb{R}}|u|^{2m+2}dx+\|h\|^{2}+\mathcal{B}_{1}]ds+2\sum\limits_{j=1}^{\infty}b_{j}\int_{0}^{t\wedge \tau_{N}}(\dot{u}+\alpha u,e_{j})d\beta_{j}(s).
\end{split}
\end{equation*}
By taking the limit (as $N\rightarrow\infty$) in the above equality, we have
\begin{equation*}
\begin{split}
\mathcal{E}_{u}(t)=&\mathcal{E}_{u}(0)+\int_{0}^{t}[-2\alpha\|u\|_{1}^{2}+2(\alpha-\gamma)\|\dot{u}+\alpha u\|^{2}-2\alpha(\alpha-\gamma)(u,\dot{u}+\alpha u)\\
&-2\alpha\int_{\mathbb{R}}|u|^{2m+2}dx+\|h\|^{2}+\mathcal{B}_{1}]ds+2\sum\limits_{j=1}^{\infty}b_{j}\int_{0}^{t}(\dot{u}+\alpha u,e_{j})d\beta_{j}(s).
\end{split}
\end{equation*}
By taking the expectation of the above inequality, we have
\begin{equation*}
\begin{split}
\mathbb{E}\mathcal{E}_{u}(t)=&\mathbb{E}\mathcal{E}_{u}(0)+\mathbb{E}\int_{0}^{t}[-2\alpha\|u\|_{1}^{2}+2(\alpha-\gamma)\|\dot{u}+\alpha u\|^{2}-2\alpha(\alpha-\gamma)(u,\dot{u}+\alpha u)\\
&-2\alpha\int_{\mathbb{R}}|u|^{2m+2}dx+\|h\|^{2}+\mathcal{B}_{1}]ds,
\end{split}
\end{equation*}
this implies that
\begin{equation}\label{44}
\begin{split}
\partial_{t}\mathbb{E}\mathcal{E}_{u}(t)
=&\mathbb{E}[-2\alpha\|u\|_{1}^{2}+2(\alpha-\gamma)\|\dot{u}+\alpha u\|^{2}-2\alpha(\alpha-\gamma)(u,\dot{u}+\alpha u)-2\alpha\int_{\mathbb{R}}|u|^{2m+2}dx]\\
&+\|h\|^{2}+\mathcal{B}_{1},
\end{split}
\end{equation}
by taking $\alpha$ small enough, we have
\begin{equation*}
\begin{split}
-2\alpha\int_{\mathbb{R}}(u^{2}+u_{x}^{2})dx+2(\alpha-\gamma)\|\dot{u}+\alpha u\|^{2}-2\alpha(\alpha-\gamma)(u,\dot{u}+\alpha u)\leq
-\frac{7}{4}\alpha|y|_{\mathcal{H}}^{2}.
\end{split}
\end{equation*}
Plugging this into \eqref{44}, we have
\begin{equation*}
\begin{split}
\partial_{t}\mathbb{E}\mathcal{E}_{u}(t)+\frac{3\alpha}{2}\mathbb{E}\mathcal{E}_{u}(t)\leq\|h\|^{2}+\mathcal{B}_{1}.
\end{split}
\end{equation*}
With the help of Gronwall inequality, we can prove \eqref{45}.
\par
This completes the proof.
\end{proof}

\begin{proof}[Proof of Proposition \ref{Keypro}]
In order to prove Proposition~\ref{Keypro}, we borrow some ideas from \cite[Lemma 2.1]{Chow1} and adapt them to our setting. We first recall the standard Friedrichs mollifier. Let $\rho_\varepsilon$ be a nonnegative, even $C^\infty$-function supported in the ball of radius $\varepsilon$ centered at the origin, with
$\int_{\mathbb R} \rho_\varepsilon(x)dx = 1.$
The mollification of a function $g$ is then defined by
$$
g^\varepsilon(x) := (\rho_\varepsilon * g)(x) = \int_{\mathbb R} \rho_\varepsilon(x - y) g(y) dy.
$$
We now apply the standard Friedrichs mollifier $\rho_\varepsilon *$ to \eqref{LWE}, thereby obtaining a regularized (mollified) version of the system
\begin{eqnarray}\label{LWES}
 \begin{array}{l}
 \left\{
 \begin{array}{llll}
\ddot{u^{\varepsilon}}+Au^{\varepsilon}+\gamma \dot{u^{\varepsilon}}=f^{\varepsilon}+\eta^{\varepsilon}
 \\u^{\varepsilon}(x,0)=u_{0}(x)
 \\\dot{u^{\varepsilon}}(x,0)=u_{1}(x)
 \end{array}
 \right.
 \end{array}
 \begin{array}{lll}
 {\rm{in}}~\mathbb{R}\times(0,+\infty),\\
 {\rm{in}}~\mathbb{R},
 \\{\rm{in}}~\mathbb{R}.
\end{array}
\end{eqnarray}
Since the mullified functions $u^{\varepsilon}$ is smooth in $x$, then the equation \eqref{LWES} is to be interpreted as a stochastic differential equation for each $x\in \mathbb{R}$. For the general case, we can use a standard limit process as in \cite[Lemma 2.1]{Chow1} to finish it. Thus, without loss of generality, we can assume that the solution $u$ to \eqref{LWE} is smooth in $x\in \mathbb{R}.$ For any fixed $x\in \mathbb{R},$ we can apply Ito formula to $u^{2}+u_{x}^{2}+(\dot{u}+\alpha u)^{2}$ and $\psi^{2} u^{2}+\psi^{2} u_{x}^{2}+\psi^{2} (\dot{u}+\alpha u)^{2}$, then, we integrate the results on the $\mathbb{R}$ to prove Proposition \ref{Keypro}.
\par
This completes the proof.
\end{proof}

\noindent \footnotesize {\bf Acknowledgements.}
\par
Peng Gao would like to thank Professor Kuksin S, Shirikyan A, Nersesyan V and Zhao M for fruitful discussions.
A part of this paper was finished when Peng Gao was visiting Southern University of Science and Technology, Peng Gao thanks for Professor Li LY's invitation and hospitality.
Peng Gao thanks Professor Kuksin S for his invitation to Universit\'{e} Paris Cit\'{e}, and for his scientific supervision on ergodicity and mixing for random dynamical systems. Peng Gao thanks the financial support of the China Scholarship Council (No. 202406620219). This work is supported by National Natural Science Foundation of China (Grant No. 12371188).
\par~~
\par
\noindent \footnotesize {\bf Data Availability.} \par
Data sharing not applicable to this article as no data sets were generated or analysed during the current study.
\par~~
\par
\noindent \footnotesize {\bf Competing Interests.} \par
The author declares that there is no conflict of interest.
\par
\end{document}